\documentclass[journal]{IEEEtran}
\ifCLASSINFOpdf
  \usepackage[pdftex]{graphicx}
\else
  \usepackage[dvips]{graphicx}
\fi

\usepackage[cmex10]{amsmath}
\usepackage{amsthm}

\usepackage{graphicx, bm, bbm, latexsym, amsfonts, epsfig, verbatim, epsfig, verbatim, enumerate, multirow, caption, graphicx, subcaption, algorithm, hyperref, url, multicol, color, latexsym, amsfonts, epsfig, verbatim, tikz, cases, diagbox,multicol,fixmath,makecell,booktabs,comment}
\hypersetup{
	colorlinks=true, 
	breaklinks=true,
	urlcolor= black, 
	linkcolor= black, 
	citecolor=red, 
	pdftitle={},
	pdfauthor={},
}

\usepackage[flushleft]{threeparttable}
\usepackage{graphicx, bbm, latexsym, amsfonts, epsfig, verbatim, amsmath, color, epsfig, verbatim, enumerate, multirow, caption, graphicx, subcaption, algorithm}
\usepackage[noend]{algpseudocode}
\usepackage{url, multicol, color, latexsym, amsfonts, epsfig, verbatim, tikz, cases, bbm}
\usepackage[top=1.0in, bottom=1.0in, left=1.0in, right=1.0in]{geometry}
\usepackage{diagbox,multicol,graphicx}
\usepackage{booktabs, caption, makecell}

\usepackage{threeparttable}
\usepackage{array}
\newcolumntype{L}[1]{>{\raggedright\let\newline\\\arraybackslash\hspace{0pt}}m{#1}}
\newcolumntype{C}[1]{>{\centering\let\newline\\\arraybackslash\hspace{0pt}}m{#1}}
\newcolumntype{R}[1]{>{\raggedleft\let\newline\\\arraybackslash\hspace{0pt}}m{#1}}

\newcommand{\mc}[1]{\mathcal{#1}}
\newcommand{\mb}[1]{\mathbb{#1}}

\renewcommand{\i}{{\bf i}}

\newcommand{\empi}{\hat{\mb{P}}}

\newtheorem{rema}{Remark}

\newtheorem{prop}{Proposition}

\makeatletter
\newcommand{\subalign}[1]{%
  \vcenter{%
    \Let@ \restore@math@cr \default@tag
    \baselineskip\fontdimen10 \scriptfont\tw@
    \advance\baselineskip\fontdimen12 \scriptfont\tw@
    \lineskip\thr@@\fontdimen8 \scriptfont\thr@@
    \lineskiplimit\lineskip
    \ialign{\hfil$\m@th\scriptstyle##$&$\m@th\scriptstyle{}##$\hfil\crcr
      #1\crcr
    }%
  }%
}
\makeatother

\allowdisplaybreaks

\begin{document}
\title{Distributionally Robust Decentralized Volt-Var Control with Network Reconfiguration
}

\author{Geunyeong Byeon~\IEEEmembership{Member,~IEEE},~Kibaek~Kim~\IEEEmembership{Senior Member,~IEEE}
\thanks{G. Byeon is with the School of Computing and Augmented Intelligence, Arizona State University, Tempe, AZ, USA (contact: geunyeong.byeon@asu.edu)}
\thanks{K. Kim is with the Mathematics and Computer Science Division, Argonne National Laboratory, Lemont, IL, USA (contact: kimk@anl.gov).}
\thanks{
This material is based upon work supported by the U.S. Department of Energy, Office of Science, under contract number DE-AC02-06CH11357.
}
}

\maketitle

\begin{abstract}  
    This paper presents a decentralized volt-var {\color{black}optimization (VVO)} and network reconfiguration strategy to address the challenges arising from the growing integration of distributed energy resources, particularly photovoltaic (PV) generation units, in active distribution networks. 
    To reconcile control measures with different time resolutions and empower local control centers to handle intermittency locally, the proposed approach leverages a two-stage distributionally robust optimization; decisions on slow-responding control measures and set points that link neighboring subnetworks are made in advance while considering all plausible distributions of uncertain PV outputs. 
    We present a decomposition algorithm with an acceleration scheme for solving the proposed model. Numerical experiments on the IEEE 123 bus distribution system are given to demonstrate its outstanding out-of-sample performance and computational efficiency, which suggests that the proposed method can effectively localize uncertainty via risk-informed proactive timely decisions.
\end{abstract}


%
\IEEEpeerreviewmaketitle

\section{Introduction}
\IEEEPARstart{D}{istributed} energy resources (DERs), particularly photovoltaic (PV) systems, are increasingly causing uncertain and intermittent influxes of electric power at the edge of distribution grids. Traditional operating practices for distribution systems, which monitor and control a small subset of network components on a slow timescale, are inadequate to handle instant voltage fluctuations and overvoltage problems caused by the volatile bidirectional power flow, which requires active controls based on optimal power flow (OPF) \cite{jin2020local}. Because of the limitations, distribution operators often limit the capacity of PV systems, known as hosting capacity \cite{ismael2019state}.
 
One of the key reliability requirements that limit PV integration is a national standard for voltage regulation, such as ANSI C84.1, which demands that voltages be maintained within safe limits across the grid. In order to regulate voltages, on-load tap changers (OLTCs) and switchable capacitor banks (CBs) have been traditionally used. OLTCs alter the voltage on the secondary winding, while CBs provide reactive power near demand nodes. Recently, network reconfiguration via line switching and reactive power support by smart inverters have emerged as promising control measures for handling the DER-related issues
\cite{zhou2021three}; the new IEEE 1547-2018 standard even stipulates the reactive power support capability of smart inverters for DER integration.

Recent studies have developed decentralized approaches for effective operations of such legacy and emerging control devices in the distribution grid (e.g., \cite{antoniadou2017distributed,li2018decentralized,li2019distributed,feng2017decentralized}).
These approaches aim to balance between centralized and distributed methods by assuming a hierarchical control system that consists of local control centers (LCCs) and a central coordinator (CC).
Specifically, the CC makes the slow timescale decisions on OLTCs, capacitors, and/or line switches and also determines some set-points or policies for coordinating the LCCs, each of which independently makes operational decisions on local smart inverters at a much higher time resolution to respond to the intermittent voltage changes. 

In order to define the LCCs' governing regions, network partitioning methods are proposed in \cite{li2018decentralized,li2019distributed}, and a componentwise decomposition is employed in \cite{feng2017decentralized}. In addition, the CC and/or LCCs are often modeled as optimization problems \cite{li2019distributed} or a machine learning architecture/reinforcement learning (RL) agents for fast control and coordination \cite{sun2021optimal,sun2021two,cao2021deep,liu2021robust}. 
However, none of the literature considers network reconfiguration, and RL-based approaches face challenges such as scalability, safety, and robustness, as discussed in \cite{chen2022reinforcement}. 

Uncertain intermittency of PV generation further complicates the distribution systems control.
A risk-aware OPF-based decentralized control was proposed in \cite{li2019distributed} based on a two-stage adaptive robust optimization (TSARO) problem that considers an interval estimate of uncertain PV outputs and demand. The solutions to TSARO can be conservative, however, as it ignores the likelihood of each outcome \cite{ben2002robust}. 
Alternatively, two-stage distributionally robust optimization (TSDRO) provides a more flexible approach for making risk-informed decisions for the situations when probability distributions of uncertain factors are hard to specify because of a limited number of samples or nonstationarity. TSDRO hedges against a worst-case probability distribution within a set of plausible probability distributions, namely, an ambiguity set.
Several TSDRO models are proposed for dispatch and/or network reconfiguration in \cite{zhou2020linear,zheng2020adaptive,zhou2021three,chen2021fast,liu2021data,zhai2022distributionally,maghami2023two}, but none of the work accounts for decentralization. 
In \cite{huang2020distributionally}, decentralized dispatch is considered in a TSDRO setting, where an ambiguity set is constructed on a finite sample space. 

\begin{figure*}[ht!]
    \centering
    \subfloat[Central control decisions made ahead of time by CC\label{1a}]{\includegraphics[width=0.45\textwidth]{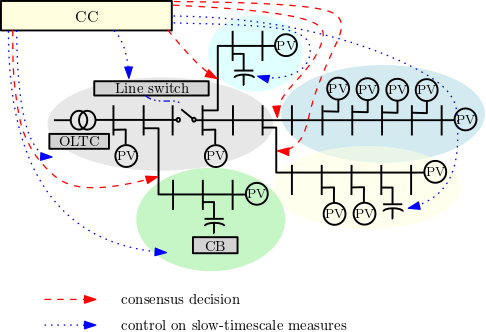}}
        \hfill
    \subfloat[Real-time local controls by LCCs \label{1b}]{\includegraphics[width=0.45\textwidth]{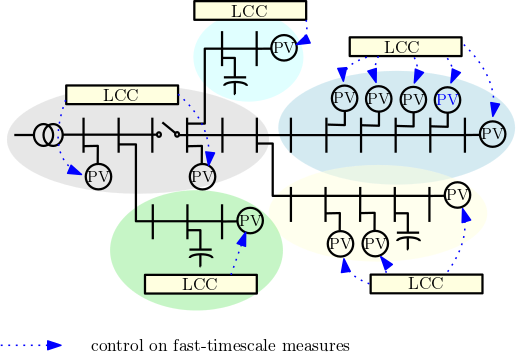}}
    \caption{Proposed decentralized control scheme} \label{fig:communication-network}
\end{figure*}

{\color{black} Driven by the need to tackle coordination challenges in the face of high uncertainty,} in this paper, we propose a TSDRO problem for the decentralized {\color{black} volt-var optimization (VVO)} and reconfiguration of distribution systems with uncertain PV outputs. The ambiguity set is
constructed with all the probability distributions that are
within a certain Wasserstein distance from a reference
distribution, namely the Wasserstein ambiguity set, to account for the uncertain PV outputs. 
Figure \ref{fig:communication-network} illustrates the proposed TSDRO approach; 
{\color{black}
\begin{itemize}
    \item The CC, in the first stage, decides on slow-responding control measures and determines the set points for variables linking neighboring subnetworks while considering the ambiguity set;
    \item Each LCC, in the second stage, operates its region by using PV inverters while meeting the set points.
\end{itemize}}

{\color{black}The coordination among OLTC, CB, line switches, and PV inverters is a result of well-considered decisions made by the CC. Specifically, each LCC autonomously manages PV inverters to regulate its respective subregion while adhering to the set points set by the CC. The CC considers the uncertainty of PV outputs and the LCCs' responses to this uncertainty when establishing set points and making slow timescale decisions. This consideration is facilitated via DRO, effectively safeguarding against the expected local system cost of LCCs under worst-case probability distributions of intermittent PV outputs from the ambiguity set.}
\noindent 
The Wasserstein ambiguity set is particularly a good choice for this application since it is robust even with a limited number of samples and can hedge against potential bias in the probabilistic forecast \cite{li2020review}.
{\color{black}To the best of our knowledge, a research gap exists regarding \emph{the application of DRO to enhance the coordination of subnetworks}. This gap holds the potential to allow reliable autonomous operations of each subnetworks so that the real-time uncertainty can be effectively absorbed within each subnetworks.}

The key contributions of this paper are threefold. First, it introduces a {\color{black}new TSDRO model} that effectively coordinates multiple LCCs and various control measures in the presence of intermittent PV generation. {\color{black}Notably, it presents an alternative formulation for reconfiguring the network as a forest.} Second, it proposes a scheme that enhances the solution approach developed in \cite{byeon2022two} for solving the TSDRO model. Unlike the heuristic solution presented in \cite{li2019distributed} for their TSARO model, the solution approach provides an exact solution to the TSDRO model. 
Third, the paper presents numerical results from a case study that demonstrate the potential advantages of the TSDRO model over two-stage robust optimization (TSRO) and the sample average approximation (SAA) of two-stage stochastic optimization (TSSO) in terms of out-of-sample performance. Specifically, the proposed method, even with a limited sample size, ensures reliable load shedding and power import in a majority of scenarios via an efficient and effective utilization of PVs.

The rest of this paper is organized as follows. Section \ref{sec:preliminaries} presents notations and preliminaries, and Section \ref{sec:centralized} formalizes the decentralized control problem. Section \ref{sec:solution} briefly reviews a solution method, and Section \ref{sec:case-study} analyzes the out-of-sample performance of the model on a test system.  Section \ref{sec:conclusion} concludes the paper.

\section{Notations and Preliminaries}\label{sec:preliminaries}
The parameters and variables are summarized in Tables \ref{table:param} and \ref{table:var}, respectively. We denote sets with calligraphic letters (e.g., $\mc N$ and $\mc E$) and use capital letters to denote matrices, unless otherwise stated.
We denote random numbers or vectors with the tilde and their realizations without the tilde. For a set $\mc A$ and a bus $i \in \mc N$, we let $\mc A(i)$ denote the subset of $\mc A$ that is associated with $i$ (e.g., $\mc K(i)$ is the set of distributed generators located at $i$). 
For an integer $n$, $[n]$ denotes a set $\{1,\cdots, n\}$. 
Sets with subscript $i$  denote its subset associated with subregion $i$; for example, $\mc S_i$ denotes the set of shunt capacitors located in subregion $i$. We assume that the uncertain PV output levels $\tilde \xi_i$ associated with subregion $i$ are independent of $\tilde \xi_j$ for $j \neq i$. 

\begin{table}[!t]\fontsize{9}{9}\selectfont
\renewcommand{\arraystretch}{1.2}
	\centering
	\caption{Parameters} \label{table:param}
	{\begin{tabular}{p{0.27\linewidth} p{0.65\linewidth}}
		\toprule
		Notation & Description\\
		\midrule
        $\mc G = (\mc N,\mc E)$ & undirected graph representing the distribution grid\\
        $\quad \mc N$ & set of buses, indexed by $j \in \{1, \cdots, |\mc N|\}$, where $1$ denotes the substation bus\\
        $\quad \mc E$ & set of lines, indexed by $\{(j,k) \in \mc N \times \mc N: j< k\}$\\
        $\mc E^u \subseteq \mc E$ & set of switchable lines\\
        $\mc K = \mc K^D \cup \mc K^{PV}$ & set of energy resources, where $\mc K^{D}$ and $\mc K^{PV}$ denote the set of dispatchable units and $k$ residential PVs, respectively \\
        $\overline g_{j}^p + {\bf i} \overline g_j^q$ & upper bound on power output of $j \in \mc K^D$\\
        $\underline g_{j}^p + {\bf i} \underline g_j^q$ & lower bound on power output of $j \in \mc K^D$\\
        $\overline p_l$ & installed capacity of $l \in \mc K^{PV}$\\
        $\mc S$ & set of shunt capacitors\\
        $\overline q_k$ & reactive power rating of $k \in \mc S$\\
        $\overline v_j, \underline v_j$ & upper and lower limits of voltage magnitude squared at $j \in \mc N$\\
        $d^p_j + {\bf i} d^q_j$ & load on $j \in \mc N$\\
        $\overline\ell_{jk}$ &  limit on current magnitude squared passing through $(j,k) \in \mc E$\\
        $r_{jk} + {\bf i} x_{jk}$ & complex impedance of $(j,k) \in \mc E$\\
        $\beta_v, \beta_{shed}, \beta_{slack}$ & weights on objective terms for substation voltage magnitude, and consensus violation, respectively\\
        \multicolumn{2}{l}{{Network decomposition}}\\
        $\mc L$ & set of LCCs\\
        $\mc{N}_i, \mc{E}_i$ & set of buses and lines of the subnetwork governed by LCC $i \in \mc L$, respectively\\
        $\mc C \subseteq \mc E$ & set of lines connecting each pair of
neighboring subnetworks\\
        $k_i$ &  number of PVs in the subregion $i \in \mc L$\\
        \bottomrule
    \end{tabular}}
\end{table}

\begin{table}[!t]\fontsize{9}{9}\selectfont
\renewcommand{\arraystretch}{1.2}
	\centering
	\caption{Variables} \label{table:var}
	{\begin{tabular}{p{0.25\linewidth} p{0.69\linewidth}}
		\toprule
		Notation & Description\\
		\midrule
		\multicolumn{2}{l}{{Continuous variables}}\\
        $v_{j}$ & voltage magnitude squared at $j \in \mc N$\\
        $p_{jk} + {\bf i} q_{jk}$& complex power flow toward $(j,k) \in \mc E$ at $j$\\
        $\ell_{jk}$ & current magnitude squared on $(j,k) \in \mc E$\\
        $g^p_{j} + {\bf i} g^q_{j}$ & power generation of $j \in \mc K$\\
        ${\bf i} q^c_{k}$ & reactive power generated by $k \in \mc S$\\
        $f_{jk}$ & artificial flow between $(j,k) \in \mc E$\\
        ${\theta}^p_{j} + {\bf i} {\theta}^q_{j}$ & amount of load shed at $j \in \mc N$\\
        $\kappa_{jk}, \iota_{jk}, o_{jk}$ & auxiliary variables for converting \eqref{eq:2nd:power-flow:power:conv} into 3D Lorentz cone constraints\\ 
        \multicolumn{2}{l}{{Binary variables}}\\
        $u_{jk}$ & 1 if $(j,k) \in \mc E^u$ is open (i.e., its end buses are not adjacent), 0 otherwise\\
        $w_{k}$ & 1 if $k\in \mc S$ is on, 0 otherwise\\
        $s_{jk}$ & adjacency of $j ,k \in \mc N$\\
        $s_{0j}$ & 1 if $j \in \mc N$ receives a nontrivial artificial flow from the dummy node, 0 otherwise \\
        \multicolumn{2}{l}{{Random variables}}\\
        $\tilde \xi_{i}=(\tilde \xi_{il})_{l\in[k_i]}$ & random output level of PVs in subregion $i$\\
        \bottomrule
    \end{tabular}}
\end{table}


In this paper we pose the control problem as a TSDRO problem that features two types of decisions: a first-stage CC's action $x$ 
that is made before the uncertain PV output levels $(\tilde \xi_i)_{i \in \mc L}$ 
are realized and a second-stage LCC's action $y_i$ 
that is made to support $x$ after observing a realization $\xi_i$ of $\tilde \xi_i$ for all $i \in \mc L$. For each $i \in \mc L$, we let 
$\mc P_i$ represent a set of plausible probability distributions of $\tilde \xi_i$, which is referred to as an \emph{ambiguity set}. The aim of TSDRO is to find $x$ that hedges against a worst-case probability distributions of $\tilde \xi_i$'s among those in $\mc P_i$'s.


\begin{table}[!t]\fontsize{9}{9}\selectfont
    \renewcommand{\arraystretch}{1.2}
        \centering
        \caption{Parameters for Wasserstein ambiguity set} \label{table:param:W}
        {\begin{tabular}{p{0.37\linewidth} p{0.56\linewidth}}
            \toprule
            Notation & Description\\
            \midrule
            $\Xi_i := [0, 1]^{k_i}$ & support of $\tilde \xi_i$ for $i \in \mc L$\\ 
            $n_i$ & sample cardinality (equiv. number of scenarios) of $\tilde \xi_i$\\
            $\zeta_{i1}, \cdots, \zeta_{in_i} \in \Xi_i$ & sample (equiv. scenarios) of $\tilde \xi_i$\\
            $P_{ij}$ & probability associated with scenario $\zeta_{ij}$\\
            $\empi_i:= \sum_{j \in [n_i]} P_{ij}\delta_{\zeta_{ij}}$ & discrete reference distribution of $\tilde \xi_i$, where $\delta_{\zeta_{ij}}$ denotes Dirac measure concentrated on $\{\zeta_{ij}\}$ \\
            $p$ & scalar from $[1,\infty]$ that defines the $l_p$-norm over $\Xi_i$ for all $i \in \mc L$\\
            $\epsilon_i$ & Wasserstein ball radius for $i \in \mc L$\\
            \bottomrule
        \end{tabular}}
    \end{table}

Specifically, for each $i \in \mc L$, we use a \emph{Wasserstein ambiguity set}, parameters of which are given in Table \ref{table:param:W}.
The Wasserstein ambiguity set is composed of probability distributions that are within $\epsilon_i$-distance of a reference distribution $\hat{\mb P}_i$ with regard to a Wasserstein metric $d(\cdot, \cdot)$, for some scalar $\epsilon_i > 0$.
Let ${\mc P_i}(\Xi_i)$ denote the collection of all probability distributions $\mb Q$ on $\Xi_i$ 
with a finite first moment. 
A \emph{Wasserstein distance} between two probability distributions $\mb Q_1$ and $\mb Q_2$ in $\mc P_i(\Xi_i)$ is defined as
$$d(\mb Q_1 ,\mb Q_2):=\inf _{\gamma \in \Gamma (\mb Q_1 ,\mb Q_2 )}\left\{\int _{\Xi_i\times \Xi_i}\|\xi_1-  \xi_2\|_p\, \mathrm{d}\gamma (\xi_1,\xi_2)\right\},$$
where 
$\Gamma (\mb Q_1 ,\mb Q_2 )$ denotes the collection of all probability distributions 
with marginals 
$\mb Q_1$ and $\mb Q_2$, respectively, and $\|\cdot\|_p$ is an $l_p$-norm in $\mb R^{k_i}$ for some $p \in [1,\infty]$.  
In this paper we {\color{black}set $p=1$ and }use a discrete distribution $\empi_i$ as the reference (e.g., the empirical distribution or a probabilistic forecast given as a histogram). 
Formally, the Wasserstein ambiguity set is defined as 
\begin{equation*}
    \mc P_i := \{\mb Q \in \mc P_i(\Xi_i): d(\mb Q, \empi_i) \le \epsilon_i\}.
\end{equation*}

\section{Distributionally Robust Reconfiguration and Decentralized Control }\label{sec:centralized}

We present a TSDRO formulation for distributionally robust decentralized
VVC and reconfiguration under uncertain PV generation. On a rolling basis, the CC solves \eqref{prob} to decide on the first-stage decision $x$ that contains the consensus decision $\mu$ for coordinating LCCs (e.g., voltages and power flows on coupling lines), together with the slow-responding control decisions $\sigma$ (e.g., OLTC operations and switches for lines and capacitor banks). The decision is made while considering the worst-case expected real-time system cost $\sup_{\mb P_i \in \mc P_i} \mb E_{\tilde\xi_i \sim \mb P_i}[Z_i(x,\tilde\xi_i)]$ for each $i \in \mc L$, where $Z_i(x,\xi_i)$ computes the real-time cost of subregion $i \in \mc L$ for given first-stage decision $x$ and a realization $\xi_i$ of $\tilde \xi_i$:
\vspace{-4mm}
\begin{subequations}
    \begin{align}
     \min_{x} \ & c^T x + \sum_{i \in \mc L} \sup_{\mathbb P_i \in \mc P_i} \mb E_{\tilde\xi_i \sim \mb P_i} [ Z_i(x, \tilde \xi_i)]\\
     \mbox{s.t.} \ & x = (\sigma, \mu),\\
     \ & \sigma = (u, v_1, w) \in \mc A,\label{prob:A}\\
     \ &  \mu = (v_k, v_l, \ell_{kl}, p_{kl}, q_{kl})_{(k,l) \in \mc C} \in \mc B,
    \end{align}\label{prob}
    \end{subequations}
\noindent where $c^T x := \beta_v v_1 +  \sum_{(i,j) \in \mathcal C} r_{ij}\ell_{ij}$ and $\mc A$ and $\mc B$ respectively denote the feasible regions of $\sigma$ and $\mu$. The detailed
definitions of $\mc A$, $\mc B$, and $Z_i$ are 
discussed in the following sections. 

{\color{black}
\begin{rema}
The inclusion of the objective term $\beta_v v_1$ aims to achieve conservation voltage reduction (CVR) to some extent when dealing with voltage-dependent loads. As discussed in \cite{farivar2011inverter, farivar2012optimal}, maximizing CVR savings is equivalent to minimizing $\sum \beta_j v_j$ over all voltage-dependent loads where $\beta_j = \frac{\alpha_j}{2} d_j^p$ for which $0 \le \alpha_j \le 2$ denotes the exponent factor of the consumption model of the load at node $j$. A more accurate CVR formulation will incorporate $\sum_{j \in \mathcal N_i} \beta_j v_j$ into the objective of each subnetwork in $Z_i(x,\tilde\xi_i)$.
\end{rema}}
\begin{rema}
    For simplicity, we model the OLTC decision, $v_1$, as a continuous variable. {\color{black}This assumption is not expected to significantly impact the case study of this paper, given that we assume one OLTC at the substation. Rounding the obtained $v_1$ to the nearest tap value, with $\Delta$ representing the change in $v_1$, and applying the same change to all subsequent $v_i$'s (i.e., $v_i' = v_i + \Delta$) preserves the satisfaction of all power flow constraints, to be defined in \eqref{eq:2nd:power-flow:power:conv}-\eqref{eq:2nd:gen-bound-dispatchable}, since only the difference in $v_i$'s matters and the topology is radial. There may be a slight change in reactive power generated by active capacitors, in \eqref{eq:2nd:shunt-capacity}, but adjusting the reactive power support of PVs, in \eqref{eq:2nd:gen-bound-PV}, may mitigate this without significantly altering the solution.

However, it is important to note that when dealing with multiple distributed voltage regulators and precise modeling of voltage-dependent loads, discrete modeling of OLTC may become imperative. The rounding procedure in such cases may lead to inapplicable solutions. 
    
    Nonetheless, we emphasize that \eqref{prob} can readily incorporate a discrete OLTC modeling as in \cite{li2019distributed}. Introducing discrete modeling of OLTCs is not anticipated to significantly impact computational difficulty. This is mainly due to the typically small number of OLTCs involved, and the set of binary variables used to model an OLTC forms a special ordered set---a structure that most solvers can efficiently handle.}
\end{rema}

\begin{figure}
    \centering
    \includegraphics[width=0.9\linewidth]{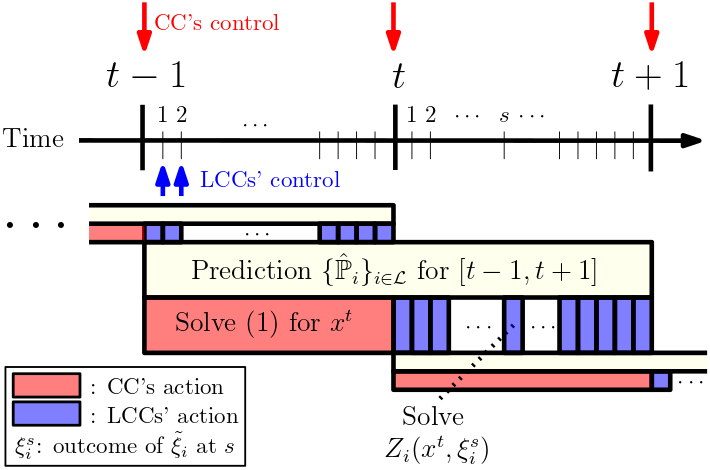}
    \caption{Rolling decision horizon of the CC}
    \label{fig:decision-rolling-horizon}
\end{figure}

\begin{rema}[Rolling decision horizon]
Problem \eqref{prob} is designed to be solved on a rolling basis. 
{\color{black}Let $t$ represent the index for the time at which the CC executes its control decision, typically on an hourly basis. LCCs, operating on a finer timescale, use the index $s$ to denote the time at which they make control decisions within two consecutive $t$'s, typically occurring every 5-15 minutes.} 
As illustrated in Figure \ref{fig:decision-rolling-horizon}, on a rolling basis, the CC is given a probabilistic forecast $\empi_i$ of $\tilde \xi_i$ for upcoming time periods{\color{black}, possibly via the workflow suggested in \cite{li2020review}. } 
With $\{\empi_i\}_{i \in \mc L}$ updated, the CC solves \eqref{prob} to obtain $x$ to be in effect in the next period, say, $x^{t}$; then, $x^{t}$ can be used in the following time period to enable autonomous operations of LCCs. Specifically, {at each $s$ within the interval $[t, t+1]$,} each LCC $i$ observes a realization $\xi^{s}_i$ of $\tilde \xi_i$ and operates its region by solving $Z_i(x^{t}, \xi^{s}_i)$ independently. 
\end{rema}
\subsection{Radiality constraints $\mc A$}\label{sec:model:first-stage}

To define $\mc A$, 
we first establish the following constraints on line switching:
\begin{subequations}
\label{eq:1st:switch}
\begin{align}
    &s_{jk} = 1, \ \forall (j,k) \in \mc E \setminus \mc E_u, \label{eq:1st:switch:line-without-switch}\\
    &s_{jk} = 1 - u_{jk}, \ \forall (j,k) \in \mc E_u, \label{eq:1st:switch:edge-on-off}
\end{align}
\end{subequations}
where $s_{jk}$ indicates whether buses $j$ and $k$ are adjacent or not. 
The vector $s$ is often referred to as a characteristic vector that represents the network configuration.


{\color{black}Radial topology is often required for the effective operation of distribution grids \cite{lei2020radiality}. As DERs become integrated, there is a demand for the network to adopt a forest configuration---a collection of connected components without cycles. To address this necessity, Lei et al. \cite{lei2020radiality} introduced a simple approach to enforce the forest constraint; it incorporates additional variables to indicate active topology and constraining them to form a subnetwork of a spanning tree. 
In this paper, we propose an alternative scheme for enforcing the forest with a smaller number of feasible solutions. 

Our idea is to enforce the spanning tree requirement on an extended graph $\overline{\mc G}=(\overline{\mc N}, \overline{\mc E})$ of $\mc G$, in which a dummy node, indexed by $0$, and a set of lines connecting the dummy node to each bus are added, i.e., $\overline{\mc N} = \mc N \cup \{0\}$ and $\overline{\mc E} = \mc E \cup \{(0,j):j \in \mc N\}$. For $j \in \mc N$, we define a binary variable $s_{0j}$, which represents the connectivity of dummy node $0$ and $j$. Thus, $\bar s := (s_{jk})_{(j,k) \in \overline{\mc E}}$ is a characteristic vector that represents a network configuration of $\overline{\mc G}$. 
Let proj$_{s}(\bar s)$ denote the projection of $\bar s$ onto the space of $s$ which is obtained by removing $s_{0j}$'s; let $\Phi(\mc G)$ denote the set of all characteristic vectors of forests of $\mc G$.

We show that the following constraint ensures the forest topology, the proof of which is in Appendix \ref{appendix:radiality}:
\begin{equation}
\bar s \in \Omega(\overline{\mc G}),\label{eq:spanning-tree}
\end{equation}
where $\Omega(\overline{\mc G})$ denotes the set of all characteristic vectors of spanning trees of $\overline{\mc G}$. 
\begin{prop}
Let $\mc S = \{\mbox{proj}_s(\bar s):\eqref{eq:spanning-tree}\}$. Then, $\mc S = \Phi(\mc G)$, i.e., $\mc S$ is the set of all characteristic vectors of forests of $\mc G$. 
\label{prop:radiality}
\end{prop}

Therefore, we can ensure a forest topology by requiring $\overline{\mc G}$ to form a spanning tree. For example, we can employ a compact formulation \cite{lavorato2011imposing} for \eqref{eq:spanning-tree}. Consider a network flow problem defined on $\overline{\mc G}$, in which the dummy node is the supply node that supplies $|\mc N|$ units of flow, and each node in $\mc N$ demands a single unit of flow:
\begin{subequations}
\label{eq:1st:network-flow}
    \begin{align}
        &0 \le f_{0j} \le |\mc N| s_{0j}, \ \forall j \in \mc N, \label{eq:1st:network-flow:1}\\
        &-|\mc N| s_{jk} \le f_{jk} \le |\mc N| s_{jk}, \ \forall (j,k) \in  \mc E, \label{eq:1st:network-flow:2}\\
        &f_{0j} + \sum_{(k,j) \in \mc E} f_{kj} - \sum_{(j,k) \in \mc E} f_{jk} = 1, \ \forall j \in \mc N,\label{eq:1st:network-flow:demand}\\
        &\sum_{j \in \mc N} f_{0j} = |\mc N|,\label{eq:1st:network-flow:supply}\\
            &\sum_{(j,k) \in \mc E} s_{jk} + \sum_{j\in \mc N} s_{0j}\le |\mc N|.\label{eq:1st:tree}
    \end{align}\label{eqs:single-commodity}
    \end{subequations}
Eqs. \eqref{eq:1st:network-flow:1}--\eqref{eq:1st:network-flow:supply} ensures the connectivity of $\overline {\mc G}$ by requiring each node to receive one artificial flow from the dummy node.
Then, Eq. \eqref{eq:1st:tree} guarantees the radial structure of $\overline{\mc G}$, which requires the number of active lines should be no more than the number of nodes minus 1. 
\begin{rema}
Eq. \eqref{eq:spanning-tree} can be modeled by any other spanning tree formulations proposed in existing literature (e.g., \cite{lei2020radiality,fan1996distribution,lavorato2011imposing,ramos2005path,jabr2012minimum,ahmadi2015mathematical}). As long as $\Omega(\overline{\mc G})$ denotes a spanning tree polytope, i.e., the convex hull of the characteristic vectors of spanning trees of $\overline{\mc G}$, the set $\mc S$ defines a forest polytope for the graph $\mc G$.
\end{rema}
\begin{rema}
A most relevant work is \cite{lei2020radiality}, where $|\mc E|$ additional binary variables, denoted by $\alpha$, are introduced to formulate $\Phi(\mc G)$ as follows: 
\begin{subequations}
\begin{align}
\mc S' = \{\alpha \in \mathbb B^{|\mc E|} : \alpha_{jk} \le s_{jk}, \ \forall (j,k) \in \mc E,\\
s \in \Omega(\mc G)\}.\label{eq:lei:spanning-tree}
\end{align}    
\end{subequations}
It is shown in \cite{lei2020radiality}, $\mc S' = \Phi(\mc G)$, and it also defines a forest polytope when $\Omega(\cdot)$ is a spanning tree polytope. 

However, in the lifted spaces of $\mc S$ and $\mc S'$ in $\bar s$ and $(\alpha,s)$, $\mc S$ has a much smaller number of feasible solutions that are all equivalent in the original space. Let $\alpha \in \mc S'$ represent a forest with $c$ connected components, and let $\mathfrak{n}_1, \cdots, \mathfrak{n}_c$ denote the number of nodes in each component. Note that the total number of spanning trees of $\mc G$ covering the forest is given by $\mathfrak{N}(\alpha):=\mathfrak{n}_1 \times \cdots \times \mathfrak{n}_c \times c^{(c-2)}$. This counts all possible combinations of selecting one node from each component and forming a tree among the chosen nodes. Therefore, there exist $\mathfrak{N}(\alpha)$ equivalent solutions in the $(\alpha,s)$-space.
In the extreme case of having $|\mc N|$ connected components each containing only one bus, the number of equivalent solutions becomes $|\mc N|^{(|\mc N| - 2)}$ in the $(\alpha,s)$-space. In contrast, our proposed approach yields a significantly lower number of feasible solutions in the $\bar s$-space. Specifically, the count is $\mathfrak{n}_1 \times \cdots \times \mathfrak{n}_c$, involving the selection of one node from each component to connect with the dummy node. This count is considerably less than $\mathfrak{n}_1 \times \cdots \times \mathfrak{n}_c \times c^{(c-2)}$ when $c \ge 2$. The number can be further reduced to 1 by introducing a small enough cost on $s_{0j}$ that implicitly orders the nodes.

Another pertinent work is \cite{arif2018dynamic}, where a specific spanning tree requirement is imposed on $\overline{\mc G}$.
\end{rema}}



To summarize, we formulate the feasible set $\mc A$ as
\begin{align*}
 \mc A := & \left\{(u, v_1, w) \in \mb B^{|\mc E_u|}\times \mb R \times \mb B^{|\mc S|}: \eqref{eq:1st:switch},\eqref{eqs:single-commodity}\right\},
\end{align*}
where $s$ and $f$ are used as auxiliary binary and continuous variables, respectively.

\subsection{Decentralized operations cost $Z_i$}\label{sec:model:second-stage}
Given $\hat x$ and $\xi_i$, each LCC $i \in \mc L$ is operated to minimize generation cost, transmission loss, and load shed by utilizing local smart PV inverters. The LCC operations problem can be formulated as a second-order cone programming (SOCP) problem that takes the SOCP relaxation of the power flow physics, as proposed in \cite{farivar2012optimal}, together with constraints on line and CB switching. 
The optimal decentralized operations cost $Z_i(\hat x, \xi_i)$  given the first-stage variables $\hat x$ and $\xi_i$ is computed by solving the following SOCP problem:
\vspace{-4mm}{\par\small\selectfont\begin{subequations}\label{2nd:opf}
 \begin{align}
   \min \ & \ \sum_{j \in \mc K_i^D} C_j (g_j^p) + \beta_{shed} \sum_{j \in \mc N_i} (\theta^p_{j}+\theta^q_{j}) +
    \sum_{(j,k) \in \mc E_i} r_{jk} \ell_{jk} \\
  \mbox{s.t.} \
    &\forall (j,k) \in \mc E_i: \ell_{jk} v_{j} \ge p_{jk}^2 + q_{jk}^2, 
    \label{eq:2nd:power-flow:power:conv}\\
    &\forall (j,k) \in \mc E_i \setminus \mc E_i^u:\nonumber\\
    &\ v_{k}
    =v_{j} - 2 (r_{jk} p_{jk} + x_{jk} q_{jk}) + (r_{jk}^2 + x_{jk}^2) \ell_{jk},\label{eq:2nd:power-flow:Ohm:conv}\\
    & \ 0 \le \ell_{jk} \le \overline \ell_{jk},
    \label{eq:2nd:bound:current:active}\\ 
    &\forall (j,k) \in \mc E^u_i:\nonumber\\
    &\ v^d_{jk} \le (\overline v_j - \underline v_k)u_{jk} + v_j-v_{k},\label{eq:2nd:power-flow:Ohm:inactive:MC1}\\
    &\ v^d_{jk} \le (\overline v_j - \underline v_k)(1-u_{jk}), \\
    &\ v^d_{jk} \ge -(\overline v_j - \underline v_k) u_{jk} +  v_j-v_{k},\\ 
    &\ v^d_{jk} \ge -(\overline v_j - \underline v_k) (1-u_{jk}),\label{eq:2nd:power-flow:Ohm:inactive:MC4}\\
    & \ v^d_{jk} = 2 (r_{jk} p_{jk} + x_{jk} q_{jk}) - (r_{jk}^2 + x_{jk}^2)  \ell_{jk},\label{eq:2nd:power-flow:Ohm:inactive}\\
    & \ 0 \le \ell_{jk} \le \overline \ell_{jk} (1-u_{jk}),\label{eq:2nd:bound:current:inactive}\\
    &\forall j \in \mc N_i:\nonumber\\
    &\ \sum_{l \in \mc K_i(j)} g^p_l + \sum_{(k,j)\in \mc E_i}  (p_{kj} - r_{kj}\ell_{kj}) \nonumber\\
    & \qquad\qquad = \sum_{(j,k)\in \mc E_i} { p_{jk}} + d^p_j - \theta_j^p,\label{eq:2nd:power-flow:balance:p}\\
    &\ \sum_{l \in \mc K_i(j)} g^q_l + q^c_j + \sum_{(k,j)\in \mc E_i}  (q_{kj} - x_{kj}\ell_{kj}) \nonumber\\
    & \qquad\qquad  = \sum_{(j,k)\in \mc E_i} { q_{jk}} +  d^q_j - \theta_j^q,\label{eq:2nd:power-flow:balance:q}\\
    & \ \underline v_j \le v_j \le \overline v_j, \label{eq:2nd:power-flow:bound:volt}\\
    & \ 0 \le \theta^p_j \le d^p_j, \ 0 \le \theta^q_j \le \max(d^q_j,0),\label{eq:2nd:bound:loadshed}\\
    &\forall j \in \mc K_i^D: \underline g^p_j \le g^p_j \le \overline g^p_j,\ \underline g^q_j \le g^q_j \le \overline g^q_j,
\label{eq:2nd:gen-bound-dispatchable}\\
    &\forall k \in \mc S_i(j), j \in \mc N_i: q_{k}^c \le \overline q_k v_j, \ q_{k}^c \le \overline q_k \overline v_{j} w_k,  \nonumber \\
    & \qquad \qquad  q_{k}^c \ge \overline q_k \overline v_{j} (w_k-1) +  \overline q_k v_j, \ q_{k}^c \ge 0,
    \label{eq:2nd:shunt-capacity}\\
    &\forall l \in \mathcal K_i^{PV}: \ g^p_l \ge 0, \ \sqrt{(g^p_l)^2 + (g^q_l)^2} \le \overline p_{l}\xi_{il},
\label{eq:2nd:gen-bound-PV}\\
& x_i = \left(v_1 \mbox{ if }1 \in \mc N_i, (u_{jk})_{(j,k) \in \mc E^u_i}, (w_k)_{k \in \mc S_i}, \right.\label{eq:2nd:xi}\\
&\left.\qquad\qquad(v_k, v_l, \ell_{kl}, p_{kl}, q_{kl})_{(k,l) \in \mc C_i}\right),\\
& x_i = A_i\hat x \label{eq:2nd:consensus},
\end{align}
\end{subequations}}
where $C_j(\cdot)$ denotes a linear cost function for each dispatchable unit $j \in \mc K_i^D$. 
Equations \eqref{eq:2nd:power-flow:power:conv}, \eqref{eq:2nd:power-flow:Ohm:conv}, \eqref{eq:2nd:power-flow:balance:p}, \eqref{eq:2nd:power-flow:balance:q}, and \eqref{eq:2nd:power-flow:bound:volt} formulate the SOCP relaxation of OPF proposed in \cite{farivar2012optimal}. In order to prevent overloading, \eqref{eq:2nd:bound:current:active} is added, which puts a limit on the magnitude of current passing through each line. In order to model line switching, \eqref{eq:2nd:power-flow:Ohm:conv}--\eqref{eq:2nd:bound:current:active} are replaced with \eqref{eq:2nd:power-flow:Ohm:inactive:MC1}--\eqref{eq:2nd:bound:current:inactive} for $(j,k) \in \mc E^u_i$; note that \eqref{eq:2nd:power-flow:Ohm:inactive:MC1}--\eqref{eq:2nd:power-flow:Ohm:inactive:MC4} are equivalent to $v^d_{jk} = (1-u_{jk})(v_j - v_{k})$ since $v^d_{jk} = v_j - v_k$ when $u_{jk} = 0$ and $v^d_{jk} = 0$ otherwise. Therefore, when the switchable line $(j,k)$ is inactive, that is, $u_{jk}=1$, $v_{jk}^d$ and $\ell_{jk}$ become zeros, and thus $p_{jk}$ and $q_{jk}$ become zeros by \eqref{eq:2nd:power-flow:power:conv}. It is only when the line is active, namely, $u_{jk} = 0$, that \eqref{eq:2nd:power-flow:power:conv} and \eqref{eq:2nd:bound:current:active} are imposed for the line. 
In addition, for the load located at bus $j \in \mc N_i$, \eqref{eq:2nd:bound:loadshed} bounds the amount of load shed by the load, and \eqref{eq:2nd:gen-bound-dispatchable} enforces generation bounds for dispatchable generators. Eq. \eqref{eq:2nd:shunt-capacity} models the reactive power support of switchable capacitors, which is equivalent to $q_{k}^c = \overline q_k v_j w_k$.
{\color{black}We let $x_i$ denote the tuple of variables associated with subnetwork $i$, which is expressed in \eqref{eq:2nd:xi}.} Eq. \eqref{eq:2nd:gen-bound-PV}  constrains the AC power outputs of each PV $l$ within its capacity, where $\xi_{il}$ represents the ratio of the DC power generated to the installed capacity, which ranges from 0 to 1, as in \cite{li2019distributed,farivar2012optimal}. Eq. \eqref{eq:2nd:consensus} represents the consensus constraint that fixes $x_i$ to be the corresponding values of $\hat x$ with a proper definition of a matrix $A_i$.

\begin{rema}
    Instead of enforcing \eqref{eq:2nd:consensus} strictly as a constraint, we penalize any of its violations with some large penalty $\beta_{slack}$ in the objective of $ Z_i(x,\tilde \xi_i)$. {\color{black}This soft enforcement ensures that $Z_i(\hat x, \xi_i)$ remains feasible for any $\hat x$, which the solution procedure can then utilize.}
    \label{rema:consensus}
    \end{rema} 
\vspace{-4mm}
\subsection{Feasible region $\mc B$ of consensus variables $\mu$}
{\color{black}The consensus decision $\mu$ includes voltages $v_k$, $v_l$, power flows $p_{kl}$, $q_{kl}$, and power loss $\ell_{kl}$ on coupling lines $(k,l) \in \mathcal C$.}
The CC decides on $\mu$ 
while making sure they obey the physical constraints:
\begin{subequations}
\begin{align}
 \mc B := & \left\{\mu= (v_k, v_l, \ell_{kl}, p_{kl}, q_{kl})_{(k,l) \in \mc C}: \right.\forall (k,l) \in \mc C,\nonumber\\ 
  & \qquad v_k, v_l, \ell_{kl}, p_{kl}, q_{kl} \in \mb R, \ \eqref{eq:2nd:power-flow:power:conv},\label{eq:opf1}\\
  & \qquad \begin{cases}
    \mbox{Eqs. } \eqref{eq:2nd:power-flow:Ohm:conv},\eqref{eq:2nd:bound:current:active} & \mbox{ if } (k,l) \notin \mc E_u,\\
    \mbox{Eqs. } \eqref{eq:2nd:power-flow:Ohm:inactive:MC1}-\eqref{eq:2nd:bound:current:inactive} & \mbox{ o.w.}
        \end{cases}\label{eq:opf2}\\
    &\qquad \left.\underline v_k \le v_k \le \overline v_k, \underline v_l \le v_l \le \overline v_l \right\},\label{eq:vb}
\end{align}
\end{subequations}
{\color{black} where Eqs. \eqref{eq:opf1} and \eqref{eq:opf2} formulate the SOCP relaxation of the power flow physics on each coupling line $(k,l) \in \mc C$ and Eq. \eqref{eq:vb} represents the limit on the voltage magnitudes at the coupling buses. }

\section{Solution Approach} \label{sec:solution}
{\color{black}
In this section, we outline a solution approach for \eqref{prob}, proposed in \cite{byeon2022two}. For notational simplicity, we algebraically express $Z_i(x, \tilde \xi_i)$ as a conic-LP program $\min_{y \in \mathcal K}\left\{q_i^T y : W_i y = h_i(x) + T_i\tilde \xi_i \right\}$, where $W_i$ and $T_i$ are some matrices, $q_i$ is a vector, $h_i(x)$ is a vector-valued affine mapping of $x$, and $\mathcal K$ is a Cartesian product of nonnegative orthants and second-order cones. The dual of $Z_i(x, \tilde \xi_i)$ is: $\max_{\pi\in\Pi_i} ( h_i(x) + T_i\tilde \xi_i )^T \pi$, where $\Pi_i := \{\pi: W_i^T \pi \preceq_{\mc K^*} q_i\}$, where $\mc K^*$ is the dual cone of $\mc K$.}

Problem \eqref{prob} can be posed equivalently as follows \cite{byeon2022two}:
\begin{subequations}
 \begin{align}
    \min_{x \in \mc X, \lambda_i \ge 0, i \in \mc L} & c^T x + \sum_{i\in \mc L} \epsilon_i \lambda_i + \sum_{j \in [n_i]} P_{ij} t_{ij}\\
    \mbox{s.t.} \ & t_{ij} \ge g_{ij}(x,\lambda_i), \forall i \in \mc L, j  \in [n_i],\label{eq:extended:epigraph}
\end{align}\label{prob:extended}
\end{subequations}
{\color{black}where $\lambda_i \in \mathbb R$ is the dual variable associated with the constraint defining the Wasserstein ambiguity set for each subnetwork $i$, $t_{ij}$ is an auxiliary scalar variable defined for each subnetwork $i$ and scenario $j$,} and $g_{ij}(x,\lambda_i)$ is the optimal objective value of 
\begin{subequations}
\begin{align}
    \max \  & h_i(x)^T \pi + b_i^T \psi_{i} - \lambda_i \|B_{ij} z_{i} -\zeta_{ij}\|\\
    \mbox{s.t.} \ &z_{i} \in \mc Z_i, \pi \in \Pi_i, \psi_{i} \in \mc{MC}(\pi,z_{i}),
\end{align}
\label{prob:sub}
\end{subequations}
{\color{black}where $B_{ij}$ is some matrix and $z_i$ is an auxiliary binary variable vector constrained by some polyhedron $\mc Z_i$. Explicit descriptions of $\mc Z_i$ and $B_{ij}$ are given in Appendix \ref{appendix:binary-representations}. The constraint $\psi_i \in \mc{MC}(\pi, z_i)$ denotes a set of linear inequalities, often referred to as McCormick envelopes \cite{hijazi2017convex}, representing $\psi_i = \pi \circ z_i$ for which $\circ$ represents componentwise multiplication of vectors.}

From the relatively complete recourse in Remark \ref{rema:consensus}, it is guaranteed that $g_{ij}(x,\lambda_i)$ is always bounded for any $x$ and $\lambda_i$. 
One can asymptotically discover the epigraph of $g_{ij}$ (i.e., \eqref{eq:extended:epigraph}) using their supporting hyperplanes, which can be found by solving \eqref{prob:sub} with ($x, \lambda_i$) fixed at some $(\hat x, \hat \lambda_i)$. This suggests Algorithm \ref{algo}, where 
($M^0$) denotes the problem obtained by relaxing \eqref{eq:extended:epigraph} from \eqref{prob:extended}. {\color{black} Finite convergence follows from the result in \cite{byeon2022two}.}
  
  \begin{algorithm}
    \caption{Algorithm for solving \eqref{prob}}\label{algo}
    \scriptsize
    \begin{algorithmic}[1]
        \State $\texttt{k} \gets 0$; $(M) \gets$($M^0$); $LB \gets -\infty$; $UB \gets \infty$
        \While{$|UB-LB| > \delta$}
          \State Solve ($M$)
          \State $v^\texttt{k}, (x^\texttt{k}, \lambda^\texttt{k}, t^\texttt{k}) \gets$ the optimal objective value and solution of ($M$) 
          \For{$i \in \mc L$}
          \For{$j  \in [n_i]$}
          \State Solve \eqref{prob:sub} with $x,\lambda_i$ fixed as $x^\texttt{k}, \lambda^\texttt{k}_i$ 
          \State $g_{ij}^{\texttt{k}} \gets$ the optimal objective value of \eqref{prob:sub}
          \State $(\pi^\texttt{k}, \psi_{i}^\texttt{k}, z_{i}^\texttt{k}) \gets$ an optimal solution of \eqref{prob:sub}; $\omega^\texttt{k} \gets B_{ij} z_{i}^\texttt{k}$ 
          \If{$t_{ij}^\texttt{k} < g_{ij}^{\texttt{k}}$}
          \State Add $t_{ij} \ge  (h_i(x) + T_i \omega^\texttt{k})^T\pi^\texttt{k} - \lambda \|\omega^\texttt{k} - \zeta_{ij}\|$ to ($M$)
          \EndIf
          \EndFor
          \EndFor
          \State $UB \gets \min\{UB, c^T x^\texttt{k} + \sum_{i\in\mc L} \epsilon_i \lambda_i^\texttt{k} + \sum_{j \in [n_i]} P_{ij} g_{ij}^{\texttt{k}}\}$
          \State $LB \gets v^\texttt{k}$; $\texttt{k}\gets \texttt{k}+1$
        \EndWhile
      \end{algorithmic}
  \end{algorithm}

To accelerate the algorithm, in line 7 of Algorithm \ref{algo}, instead of solving \eqref{prob:sub} from scratch, we first solve it with $B_{ij}z_i$ fixed at $\zeta_{ij}$ and then resolve it with $B_{ij}z_i$ fixed at $0$ to alternatively add a suboptimal cut if it cuts off the current candidate solution. 

{\color{black}
\begin{rema}
Algorithm 1 is adaptable to decentralized implementation. Specifically, Lines 6-11 within Algorithm 1 can be independently executed by each LCC in a decentralized manner, given the solution of (M) provided by the CC. Subsequently, the CC aggregates cut information from the LCCs and proceeds to solve (M). This approach alleviates the need for the CC to hold complete data, and concurrently distributes the computational burden.
\end{rema}
}

\section{Numerical Results}\label{sec:case-study}
In our numerical experiments we observed that SOCP solvers suffer from numerical instability and result in inconsistent solutions. Therefore, throughout the experiments reported in this paper, we use a linear asymptotic relaxation of the second-stage SOCP problem $Z_i(x, \xi_i)$, denoted by $\hat Z_i(x, \xi_i)$, for the experiment. Note, however, that the implications derived from the LP relaxation should align with the SOCP problem as the relaxation {\color{black} error is bounded by 1\% for each SOCP constraint}. The details of the relaxation are given in Appendix \ref{appendix:linear-relaxation}. 

In this section we assess the performance of the proposed TSDRO model {\color{black}\eqref{prob} equipped with two distinct $\epsilon$-selection criteria, denoted by \texttt{opt} and \texttt{hm}. We compare this performance with that of TSRO (denoted by \texttt{ro}) and TSSO (denoted by \texttt{saa}) models as outlined below:}
\begin{itemize}
    \item \texttt{opt}: \eqref{prob} with an optimal choice of $\epsilon$;
    \item \texttt{hm}: \eqref{prob} with $\epsilon$ chosen via a holdout method, in which we spare 20\% of the data for validation;
    \item \texttt{saa}: $\min_{x \in \mc X} \ c^T x + \sum_{i \in \mc L} \mb E_{\empi_i} [\hat Z_i(x, \tilde \xi_i)]$; and
    \item \texttt{ro}: $\min_{x \in \mc X} \ c^T x + \sum_{i \in \mc L} \max_{\xi_i \in \Xi_i} \hat Z_i(x,\xi_i)$.
\end{itemize}
{\color{black}TSSO excels in scenarios where a large number of empirical observations is available and captures the underlying uncertainty, whereas TSRO is useful in scenarios with limited data, but often makes (overly) conservative decisions. TSDRO achieves a balance between these two models with a greater out-of-sample performance. More detailed implications of these models are given in Appendix \ref{appendix:model-comparison}.}
\begin{figure}[t!]
    \centering
    \includegraphics[width=0.8\linewidth]{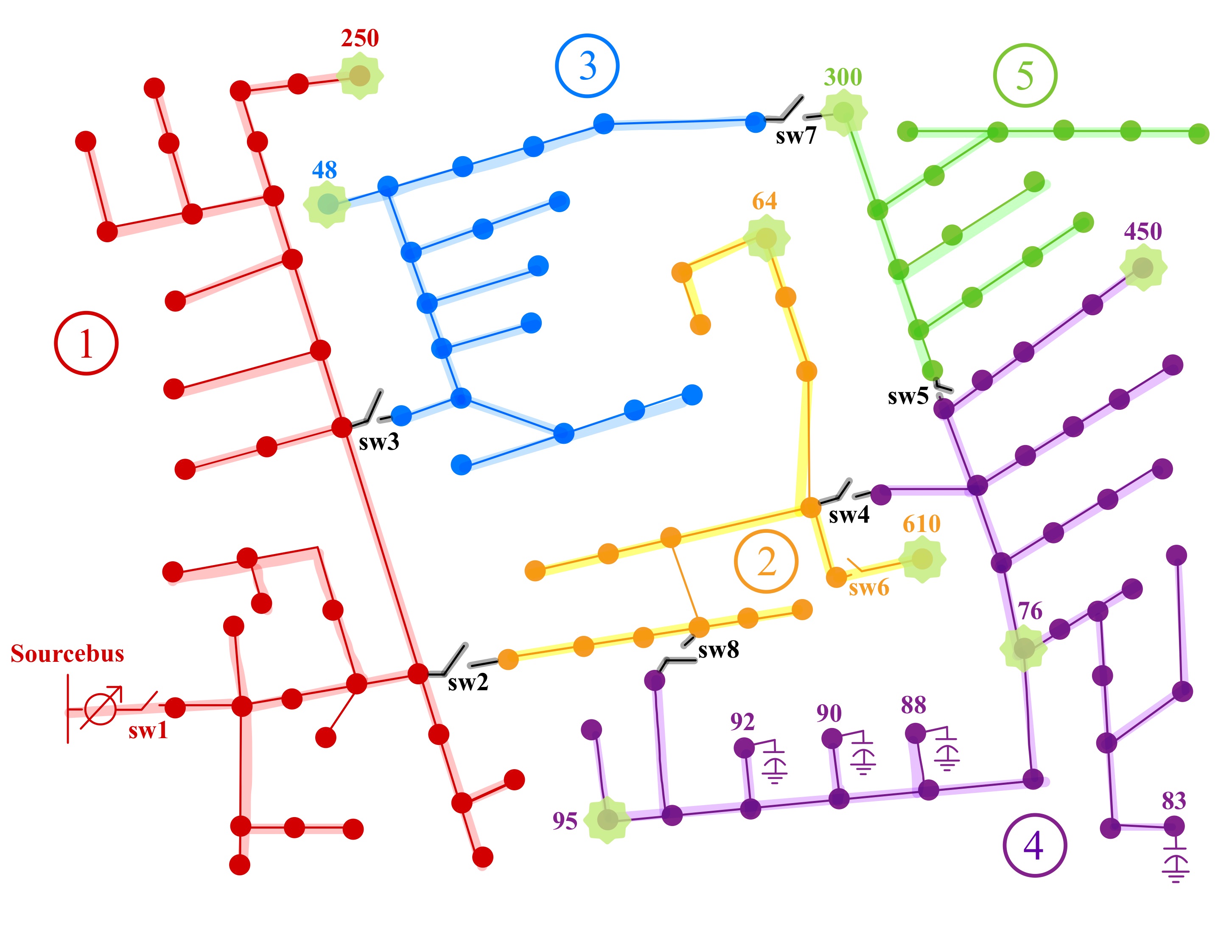}
    \caption{Modified IEEE 123 bus system with five subregions (colored subnetworks), eight PVs (green polygons), and eight switches (sw1,sw2,$\dots$,sw8).}\label{fig:test-system}
\vspace{-3mm}
\end{figure}
We use a modified IEEE 123 bus system with 8 PVs with a maximum capacity of 0.05 p.u, as illustrated in Figure \ref{fig:test-system}; see Appendix \ref{appendix:test-system} for more details.

The second stage aims to minimize the total import from the transmission grid, that is, $C_i(g_i^p) = g_i^p$ for $i \in \mc K^D$. 
The objective weights are set as $\beta_v = 0.01, \beta_{shed} = 10^2, \beta_{slack} = 10^5$ to prioritize consensus. We set the upper bound $\bar \pi$ of the dual variable associated with \eqref{eq:second-stage-constraint:xi} as $10^6$. We use $l_1$-norm to define the Wasserstein ball (i.e., $p=1$), and thus the algorithm explained in Section \ref{sec:solution} is exact. We assume 5 local control centers coupled via $\mathcal C=\{$sw2, sw3, sw4, sw5, sw7, sw8$\}$. The subregions have total loads of around $0.044 +\i 0.022, 0.032+\i0.017, 0.044+\i0.027, 0.064+\i0.035$, and $0.018+\i0.009$ p.u., respectively.

    For each subregion $i$, we assume that the true distribution of $\tilde \xi_i$ for the next time periods, denoted by $\mb Q_i$, is a truncated multivariate normal distribution \cite{doubleday2020probabilistic} with a mean of 0.8 and a covariance matrix having 0.1's on its diagonal and 0.001's on its off-diagonal entries so that PVs in the same region are slightly correlated. For the numerical experiments, the Wasserstein radii $\{\epsilon_i\}_{i \in \mc L}$ and the number of available scenarios $\{n_i\}_{i \in \mc L}$ are varied homogeneously across the subregions, so we let $\epsilon$ and $n$ denote the common radius and cardinality. For each $n \in \{5, 10, \cdots, 30\}$, we generate 50 training datasets, each of which consists of $n$ potential scenarios of $(\tilde \xi_1, \cdots, \tilde \xi_{|\mc L|})$ sampled from $\mb Q_1 \times \cdots \times \mb Q_{|\mc L|}$. Each training dataset is used to construct $\{\empi_i\}_{i \in \mc L}$ for simulation. For testing, we also sample $10^3$ out-of-sample scenarios from $\mb Q_1 \times \cdots \times \mb Q_{|\mc L|}$. The data files are available in json format on \url{https://github.com/gbyeon/DROControl-dataset.git}. 
    
    All experiments were executed on a Dell PowerEdge R650 server, with 56 cores and 512GB of RAM. The implementation is in Julia and uses IBM CPLEX 12.10. To accelerate Algorithm \ref{algo}, we add the cuts in a lazy manner to ($M$) during its branch-and-bound process using a callback function and parallelize the cut generation procedure in Lines 6-13.


\subsection{Expected total system cost}
We first compare the expected total system cost of those four models, estimated based on the $10^3$ testing scenarios, over 50 independent simulation runs. To be specific, let $\hat x$ denote a solution obtained by one of the models. Then its expected system cost is computed by $c^T \hat x + \sum_{i \in \mc L} \sum_{j=1}^{10^3} \frac{1}{10^3}\hat Z_i(\hat x, \hat\zeta_{ij})$, where $\{\hat\zeta_{ij}\}_{j=1,\cdots, 10^3}$ denotes the set of testing scenarios of subnetwork $i$.
    \begin{figure*}
    \centering
    \subfloat[$n=5$\label{1a}]{\includegraphics[width=0.33\textwidth]{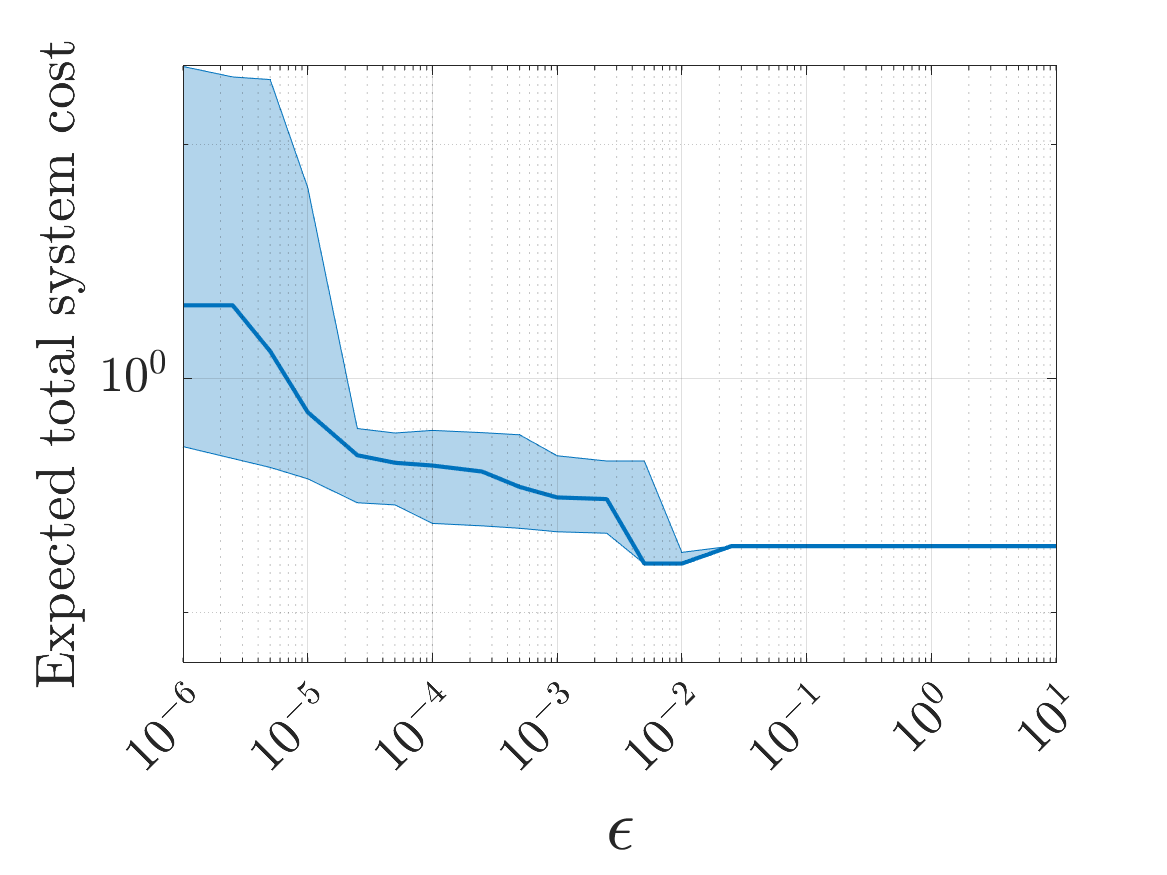}\vspace{-2mm}}
    \hfill
    \subfloat[$n=15$\label{1b}]{\includegraphics[width=0.33\textwidth]{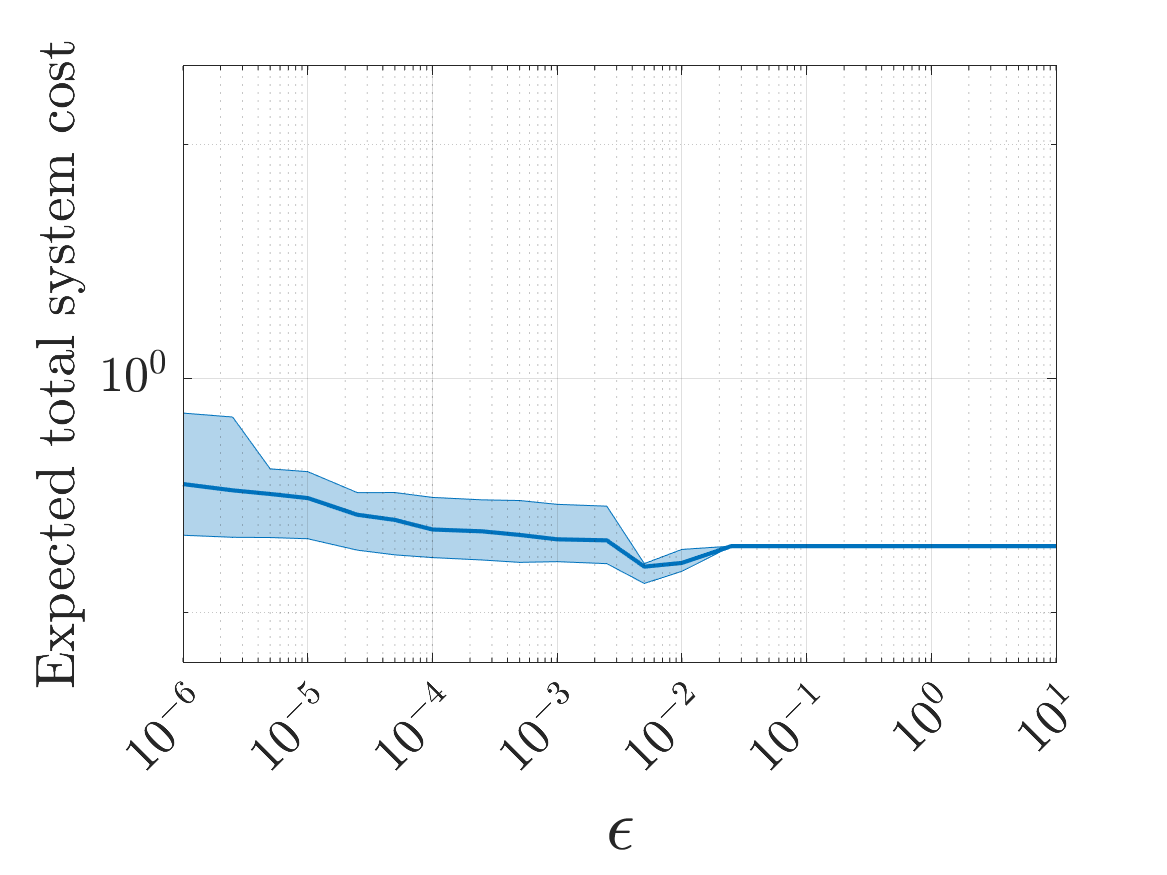}\vspace{-2mm}}
    \hfill
    \subfloat[$n=30$\label{1c}]{\includegraphics[width=0.33\textwidth]{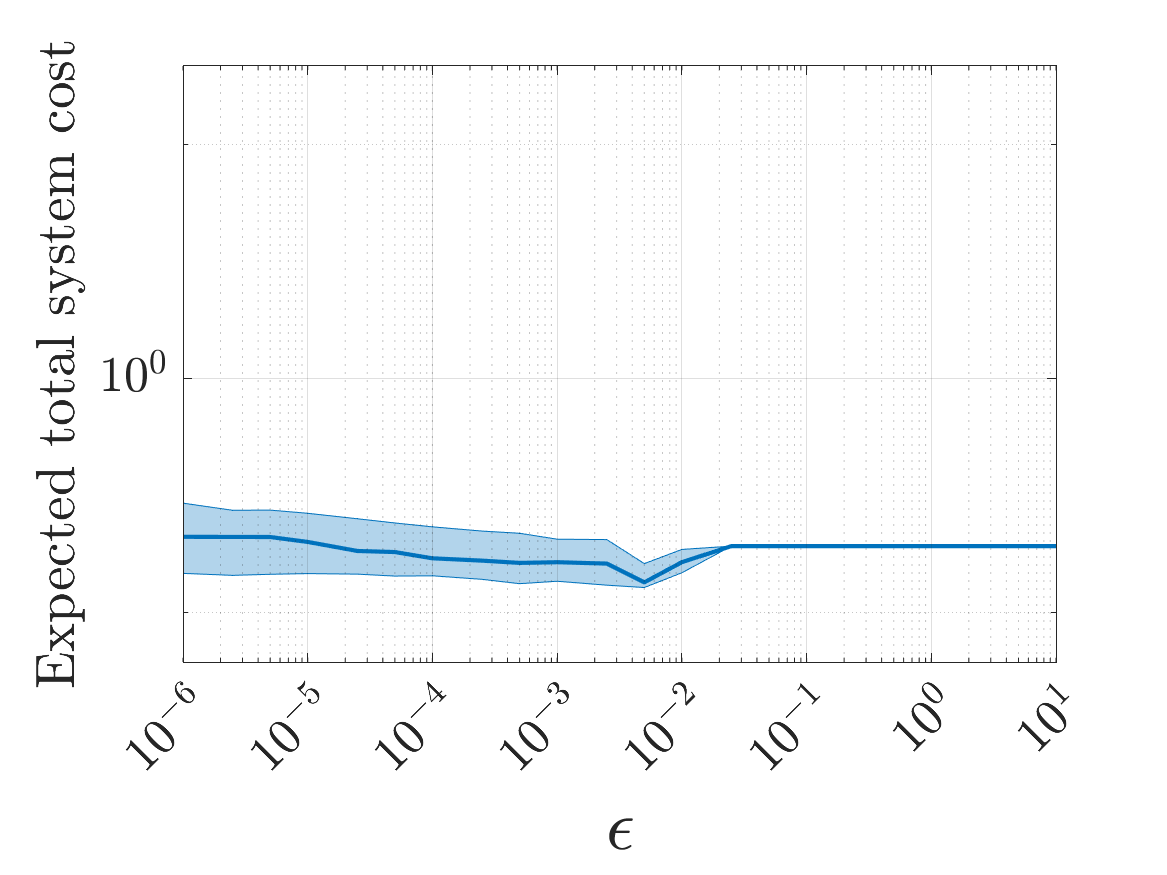}\vspace{-2mm}}
    \caption{Expected total system cost approximated by the $10^3$ testing scenarios over 50 simulations}
\label{fig:OOS_reliability}
\vspace{-4mm}
\end{figure*}
To first compare \texttt{opt}, \texttt{saa}, and \texttt{ro},
Figure \ref{fig:OOS_reliability} visualizes the 20th and 80th percentiles (shaded areas) and the means (solid lines with markers) of the expected total system cost of the solution of \eqref{prob} with respect to varying $\epsilon$, denoted by $\tilde J_{n}(\epsilon)$.
Note that as $\epsilon \rightarrow 0$, the solution of \eqref{prob} becomes the solution of \texttt{saa}; and when $\epsilon \ge 10$, the solution of \eqref{prob} becomes the solution of \texttt{ro} since the Wasserstein distance is bounded from above by the distance between two extreme scenarios $\|\boldsymbol 1 - \boldsymbol 0\|_1$, where $\boldsymbol 1$ and $\boldsymbol 0$ are vectors of ones and zeros with the dimensionality of $\max_{i \in \mc L}k_i=3$. 
In Figure \ref{fig:OOS_reliability} we observe that the total cost improves up to a certain point of $\epsilon$ and then declines for all of the simulations. This empirically shows the superiority of \texttt{opt} over \texttt{saa} and \texttt{ro}, as \texttt{opt} chooses $\tilde\epsilon$ that minimizes $\tilde J(\epsilon)$.

\begin{figure}[t!]
    \centering
    \includegraphics[width=0.7\linewidth]{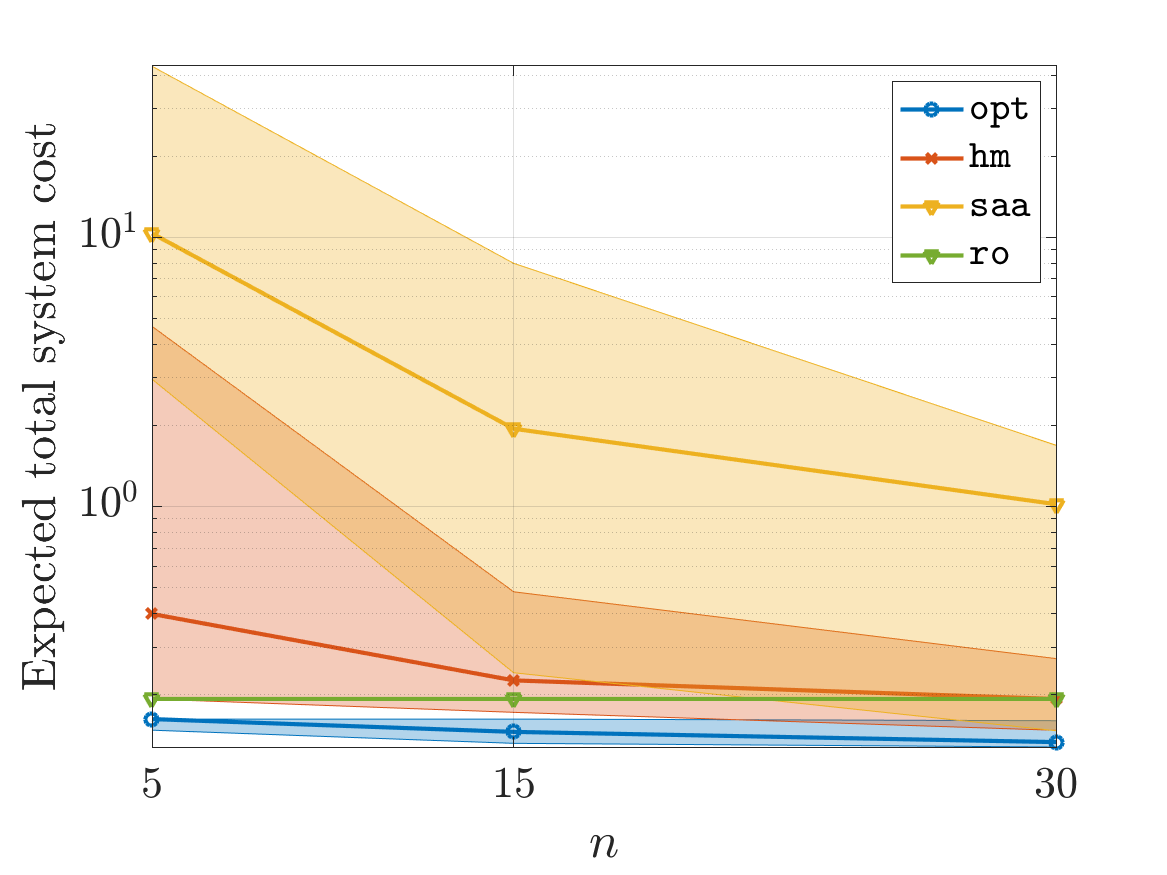}
    \caption{Expected total system cost as a function of $n$}
\label{fig:OOS_hm_opt}
\vspace{-2mm}
\end{figure}

To compare all of those four models, Figure \ref{fig:OOS_hm_opt} plots the estimated total cost with respect to the training sample size $n$. The figure highlights the sample efficiency of the DRO method. As expected, the total costs of \texttt{saa} and \texttt{hm} highly depend on the sample size $n$; one can see  from the figure that the cost improves as $n$ increases. However, since \texttt{hm} considers ambiguity in the sample, despite the suboptimal choice of $\epsilon$, its total cost is much lower than that of \texttt{saa}. On the other hand, \texttt{ro} remains the same since it does not take advantage of samples. Note that \texttt{opt} gives the minimum expected cost, even with a limited sample size. 
As expected, the performance of those three data-driven solutions \texttt{opt}, \texttt{ro}, and \texttt{saa} becomes similar as $n$ increases. 

{\color{black}
\subsection{Choice of $\epsilon$}\label{sec:eps}

 \begin{figure}
    \centering
    \includegraphics[width=0.7\linewidth]{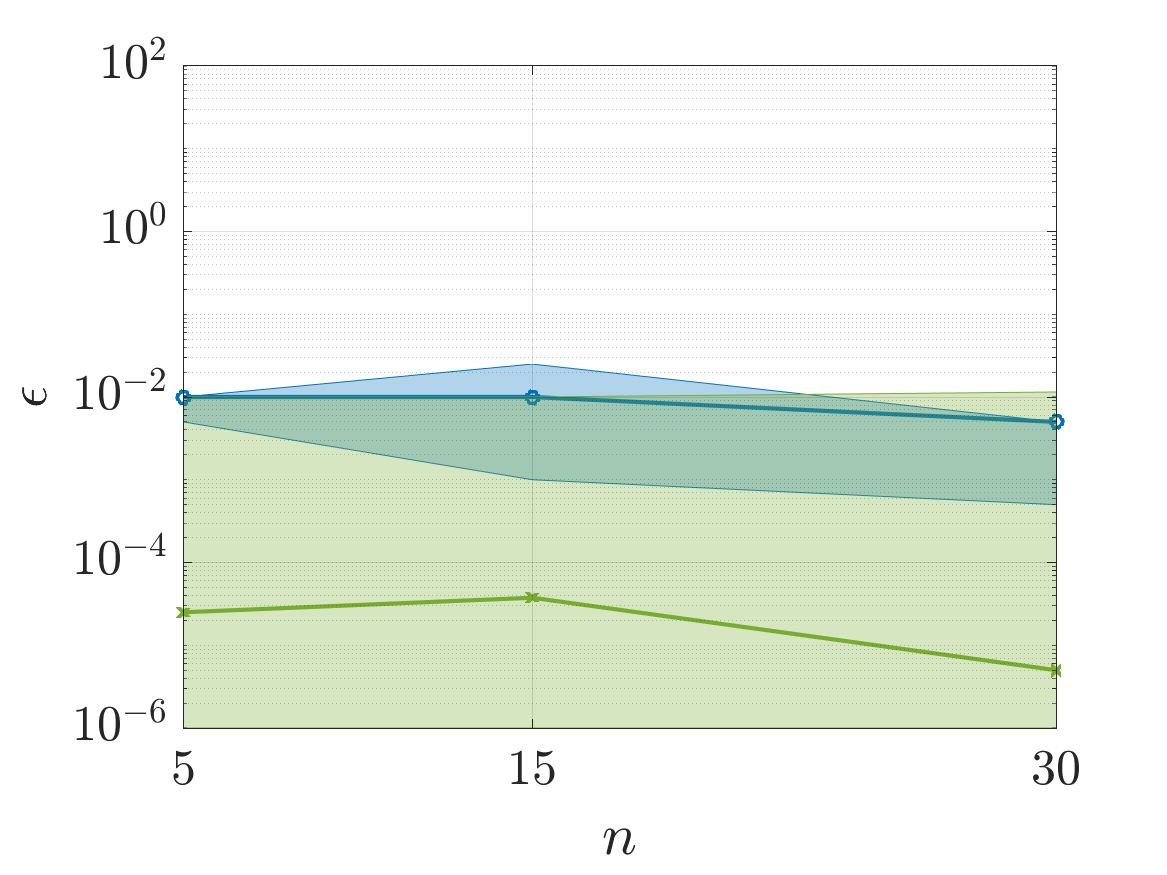}
    \caption{Choice of $\epsilon$}
    \label{fig:eps}
\end{figure}
For each cardinality choice $n$, \texttt{opt} and \texttt{hm} may use different $\epsilon$ per simulation, as explained in Appendix \ref{appendix:model-comparison}; Figure \ref{fig:eps} displays the 20th and 80th percentiles (shaded areas) and the medians (solid lines) of $\epsilon$ selected by the two methods. Figure \ref{fig:eps} indicates that the value of $\epsilon$ chosen for \texttt{opt} decreases as $n$ increases, that is, as $\{\empi_i\}$ gets more trustworthy, while that of \texttt{hm} varies significantly because of the limited amount of validation scenarios. It is worth noting that the variance in the optimal $\epsilon$ values is small for each cardinality choice. This observation suggests that learning the optimal $\epsilon$ value from experience may be promising.
}

\subsection{Impacts on load shedding}\label{sec:result:loadshed}
To see the out-of-sample performance of the four models on important system statistics, we first compare the expected networkwide load shed. 
Figure \ref{fig:loadshed_by_methods} visualizes the amount of networkwide load shed for each of the $10^3$ testing scenarios. For each scenario, the figure plots a red dot for the median and grey bars for the 10th and 90th percentiles of the results obtained over the 50 simulation runs. As Figure \ref{a} indicates, with the limited sample size \texttt{saa} can result in significant load shedding for many scenarios, and the results vary a lot per simulation, indicating its vulnerability to the limited sample size. With the suboptimal choice of $\epsilon$, the results of \texttt{hm} also fluctuate per simulation, but the red dots corresponding to the median values are much more concentrated in the lower region compared with those of \texttt{saa}. On the other hand, with the optimal choice of $\epsilon$, the results of \texttt{opt} do not vary much over different simulation runs; and, under many scenarios, there is no load shed. Lastly, \texttt{ro} has a zero load shedding for all of the testing scenarios. Figure \ref{b} shows how the results of the data-driven methods, \texttt{saa}, \texttt{hm}, \texttt{opt}, improve as we have more data. \texttt{opt} achieves zero load shedding for many of the testing scenarios, and \texttt{hm} and \texttt{saa} get much closer to \texttt{opt}. Since \texttt{ro} is not affected by the samples, it remains to have zero load shed for $n=30$.

\vspace{-2mm}
\subsection{Impacts on power imports and PV utilization}
\begin{figure*}[t!]
    \centering
    \subfloat[$n=5$\label{pva}]{\includegraphics[width=0.33\textwidth]{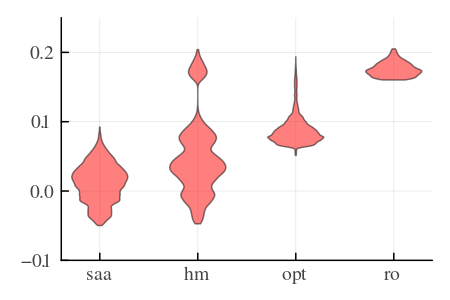}}
    \subfloat[$n=30$\label{pvb}]{\includegraphics[width=0.33\textwidth]{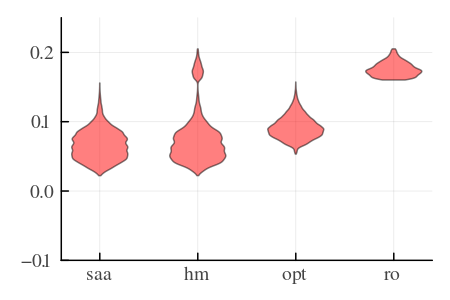}}
    \caption{Amount of real power imported at the substation}
\label{fig:pgimport}
\vspace{-2mm}
\end{figure*}
To show how much efficiency is sacrificed to prevent load shedding in those four models, we analyze PV utilization with respect to the amount of power imported from the main grid, that of utilized power outputs of PVs, and the power factor of PV outputs. Figure \ref{fig:pgimport} displays the violin plots showing the distribution of the amount of real power imported at the substation over the testing scenarios averaged for the 50 simulations. The thickness of each graph represents the density of the corresponding value. Again \texttt{ro} shows the same behavior independently of $n$, since it does not use sample data; also, its import amount is almost the total real power demand, which is 0.202. This implies that \texttt{ro} uses PV power outputs only for reactive power support, which explains why it has zero load shedding for all the scenarios. On the other hand, for $n=5$, \texttt{saa} and \texttt{hm} are often too realistic and import no more than 50\% of the total demand most of the time or even export some of the power to the main grid, resulting in a high load shed in many testing scenarios. However, \texttt{opt} imports a reliable amount that mostly covers 50\% of the demand while maintaining the level of load shed significantly lower than those of \texttt{saa} and \texttt{hm}. For $n=30$, the three data-driven methods,  \texttt{saa}, \texttt{hm}, and \texttt{opt}, become similar.

 Figures \ref{fig:util_N5} and \ref{fig:pf_N5} analyze how each method utilizes PVs; the former plots the distribution of PV utilization and the latter depicts that of the power factor. As shown in Figure \ref{fig:util_N5}, \texttt{saa} uses PVs the most while \texttt{ro} does not utilize PVs at buses 48 and 300 at all. Moreover, \texttt{ro} use PVs mostly only for reactive power support, as indicated by Figure \ref{fig:pf_N5}. On the other hand, \texttt{saa} utilizes PVs mostly for real power generation as it is too optimistic, while \texttt{hm} and \texttt{opt} do so in a more balanced way.




\vspace{-4mm}

\subsection{Computation times}
As the number of subproblems increases with the same cardinality $n$ grows, we use a varying number of cores per different $n$; 14 cores for $n=5$, 28 cores for $n=15$, 56 cores for $n=30$. Table \ref{table:comp-time} summarizes the computation time in seconds for \texttt{opt} for the 50 simulation runs. Note that the solution time can be further reduced by utilizing more cores. 
\begin{table}[!ht]\fontsize{9}{9}\selectfont
\renewcommand{\arraystretch}{1.2}
	\centering
	\caption{Computation times (sec.) for \texttt{opt}} \label{table:comp-time}
	{\begin{tabular}{R{0.1\linewidth} R{0.2\linewidth}R{0.2\linewidth}R{0.2\linewidth}}
		\toprule
		$n$ & 10th & median & 90th\\
		\midrule
  5 & 197.05 & 446.88 & 1338.50\\
  15 & 122.88& 417.89& 1212.63\\
  30 & 135.44& 285.32& 750.97\\
        \bottomrule
    \end{tabular}}
\end{table}
\vspace{-2mm}

\section{Conclusions}\label{sec:conclusion}
We proposed a decentralized control method for active distribution networks that effectively coordinates local control centers and control measures on varying timescales. This coordination is achieved by determining risk-aware here-and-now decisions centrally, which serve as set points for coupling variables and slow-responding control decisions. The here-and-now decisions are informed by a set of plausible distributions of uncertain PV outputs that are close enough to a reference distribution, such as a probabilistic forecast, in a Wasserstein distance sense. Numerical studies on the IEEE 123 bus system demonstrate the outstanding out-of-sample performance of the proposed approach; the proposed method maintained reliable load shed under a majority of testing scenarios while making the most of PVs in a balanced manner even with a limited sample size. It is demonstrated that the proposed DRO model has the potential to achieve (1) a great sample efficiency; (2) a high PV utilization level while avoiding load shed; (3) proactive coordination of LCCs and slow-responding control measures, and (4)
a successful uncertainty localization.
    This suggests that the proposed approach can be strengthened by incorporating more flexible units such as storage systems.
Future research efforts will be devoted to {\color{black}enhancing scalability for larger systems, such as the IEEE 8500-bus system, involving acceleration schemes like advanced Benders cuts \cite{fischetti2010note} and in-out approaches \cite{fischetti2017redesigning}. 
Additionally, we plan to leverage machine learning models to expedite solution procedures, particularly by predicting scenarios more likely to significantly impact the current iterate, rather than testing each sample point $\zeta_{ij}$ and the worst-case scenario.
Furthermore, our research will extend to modeling} unbalanced multiphase systems, such as the one proposed in \cite{byeon2022linear}, as well as incorporating more control devices such as storage systems. Additionally, exploring the impact of heterogeneous Wasserstein radii over subregions may be worthwhile. 

\section*{Acknowledgment}
The authors thank Dr. Dae-Hyun Choi for their invaluable comments on this paper.

\bibliographystyle{IEEEtran}
\bibliography{ref.bib}  
\nopagebreak
\begin{appendices} 
    \section{Proof of Proposition \ref{prop:radiality}}\label{appendix:radiality}


{\color{black}Suppose $s$ represents a spanning forest of $\mc G$ and let $c$ be the number of connected components in the forest. For each component $i=1,\cdots,c$, pick one node and connect it with the dummy node. Then, the resultant network will form a spanning tree of the extended graph $\mc G$. Suppose $s'$ represents a spanning tree of $\overline{\mc G}$. Then $s = proj_s(s')$ will form a subgraph of the spanning tree, which is a spanning forest \cite{lei2020radiality}.\qed}

{\color{black}
\section{Descriptions of $\mc Z_i$ and $B_{ij}$}\label{appendix:binary-representations}

\begin{multline*}\mc Z_i=\{z_i = (z_{i}^{10}, z_{i}^{11}, \cdots, z_{i}^{k_i0}, z_{i}^{k_i1}) \in \mb B^{2k_i}: \\
    z_{i}^{l0} + z_{i}^{l1} \le 1,\forall l \in [k_i]\},
\end{multline*} and 
\[B_{ij} = \left[\begin{array}{ccccccccccccccccc}
    \zeta_{ij1} & 1 \\
    & & \zeta_{ij2} & 1 \\
    & & & & \ddots\\
    & & & & & \zeta_{ijk_i} & 1
\end{array}\right].\]}




\section{A Linear Relaxation $\hat Z_i(\hat x,\xi_i)$ of $Z_i(\hat x,\xi_i)$}\label{appendix:linear-relaxation}
\begin{figure}[t!]
    \centering
    \includegraphics[width=0.85\linewidth]{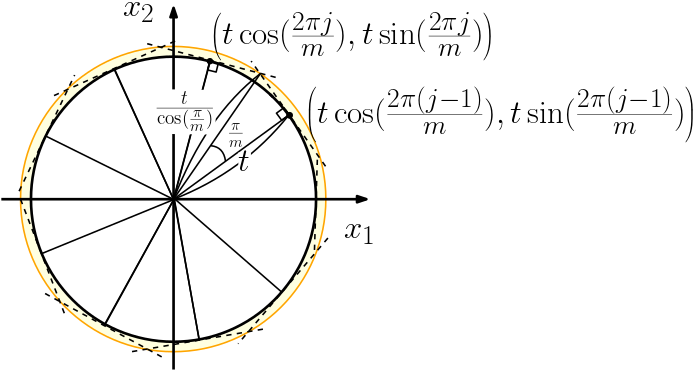}
    \caption{Horizontal slice of the linear relaxation of the 3D Lorentz cone at arbitrary $t$}\label{fig:linear-approximation}
\end{figure}

In the problem of $Z_i(\hat x,\xi_i)$, nonlinearity appears in \eqref{eq:2nd:power-flow:power:conv} and \eqref{eq:2nd:gen-bound-PV}. Note that each of these constraints can be posed as a collection of 3-dimensional Lorentz cones $\{(x,t) \in \mathbb R^2 \times \mathbb R : \sqrt{x_1^2+x_2^2} \le t \}$ and linear constraints with some auxiliary variables as follows:
\begin{subequations}
    \begin{align}
        \eqref{eq:2nd:power-flow:power:conv} \Leftrightarrow \ & \kappa_{jk} = \frac{\ell_{jk}+v_{j}}{2}, \ o_{jk} = \frac{\ell_{jk}-v_{j}}{2} \label{eq:lorentz:linear}\\
        & \iota_{jk} \ge \sqrt{q_{jk}^2 + o_{jk}^2}, \kappa_{jk} \ge \sqrt{p_{jk}^2 + \iota_{jk}^2}. \label{eq:lorentz}\\
        \eqref{eq:2nd:gen-bound-PV} \Leftrightarrow \ &  \kappa_l \le \overline p_{l}\xi_{il}, \label{eq:second-stage-constraint:xi}\\
        & \kappa_{l} \ge \sqrt{(g^p_{l})^2 + (g^q_{l})^2}.
    \end{align}
\end{subequations}

We relax each of the 3-dimensional Lorentz cone constraints with a set of $m$ linear constraints using its supporting hyperplanes at $(\hat x_1, \hat x_2, \hat t) = \left(\hat t\cos(\frac{2\pi j}{m}), \hat t\sin(\frac{2\pi j}{m}), \hat t \right)$ for $j\in [m]$ for some $m \ge 4$:
\begin{equation}
    \cos\left(2\pi j/m\right) x_1 + \sin\left(2\pi j/m\right) x_2 \le t, \ \forall j \in [m].
    \label{eq:2nd:Lorentz:linear}
\end{equation}
Figure \ref{fig:linear-approximation} illustrates the linear relaxation and indicates that the relaxation is confined in the following relaxed Lorentz cone constraint:
$$(1+\epsilon) t \ge \sqrt{x_1^2+ x_2^2},$$
where $\epsilon = \frac{1}{\cos(\frac{\pi}{m})} - 1$. As $m$ goes to infinity, $\epsilon$ approaches zero, becoming tight to the original constraint. 
\begin{rema}
    Ben-Tal and Nemirovski \cite{ben2001polyhedral} proposed a linear relaxation that approximates the 3D Lorentz cone constraint with a polynomially increasing number of additional variables and constraints. It is shown that with $2(v+1)$  additional variables and $3(v+2)$ constraints for some integer $v$, $\epsilon = \frac{1}{\cos(\frac{\pi}{2^{v+1}})}-1$ relaxation is achieved. Although the growth in problem size is at a much slower pace in the relaxation proposed in \cite{ben2001polyhedral}, \eqref{eq:2nd:Lorentz:linear} requires fewer variables and constraints up to $\epsilon = 1\%$ and is simpler; thus, we use \eqref{eq:2nd:Lorentz:linear} with $\epsilon = 1\%$ (i.e., $m = 23$). One can use the linear relaxation in \cite{ben2001polyhedral} if higher accuracy is desired.
\end{rema}

The relaxation for \eqref{eq:2nd:power-flow:power:conv} no longer guarantees $p_{jk}$ and $q_{jk}$ to be zero when $\ell_{jk} = 0$, that is, when $u_{jk} = 1$. Therefore, we add the following constraints for $(j,k) \in \mc E^u_i$:
\begin{subequations}
    \begin{align}
        -\overline s_{jk} \left(1-u_{jk}\right)\le p_{jk} \le  \overline s_{jk} \left(1-u_{jk}\right),\\
        -\overline s_{jk} \left(1-u_{jk}\right)\le q_{jk} \le  \overline s_{jk} \left(1-u_{jk}\right),
    \end{align}    \label{eq:2nd:linear:switch}
\end{subequations}
where $\overline s_{jk}$ is the apparent power limit on line $(j,k)$.

{\color{black}\section{Models for decision-making under uncertainty}\label{appendix:model-comparison}
In this paper, we compare four decision-making models under uncertainty: distributionally robust optimization (\texttt{opt}, \texttt{hm}), robust optimization (\texttt{ro}), and stochastic optimization (\texttt{saa}). \texttt{saa} finds the here-and-now decision $x$ that minimizes the total expected system cost approximated by the empirical distribution $(\empi_i)$, i.e., a minimizer of $\min_{x \in \mc X} \ c^T x + \sum_{i \in \mc L} \mb E_{\empi_i} [\hat Z_i(x, \tilde \xi_i)]$. Stochastic optimization performs well when $\empi_i$ represents the true distribution $\mb Q$ well. However, if the representation is insufficient, particularly when the sample size $n$ is limited, it tends to generate overly narrow decision $x$ that performs badly in scenarios not covered by the training sample.

In contrast, \texttt{ro} determines $x$ by hedging against the worst-case outcome of the uncertain factor, i.e., a minimizer of $\min_{x \in \mc X} \ c^T x + \sum_{i \in \mc L} \max_{\xi_i \in \Xi_i} \hat Z_i(x,\xi_i)$. Consequently, \texttt{ro} tends to be too conservative as it disregards the likelihood of outcomes, focusing solely on the worst-case scenario that may be unlikely to occur.

Model \eqref{prob} achieves a balance between \texttt{saa} and \texttt{ro} by considering an ambiguity in $\empi_i$ and considering a set of plausible distributions proximate to $\empi_i$. It hedges against the worst-case distribution within this ambiguity set. The parameter $\epsilon$ controls the range of distributions to consider. Depending on its value, the model closely resembles the solution of \texttt{saa} (when $\epsilon=0$, the ambiguity set comprises only $\empi_i$, rendering \eqref{prob} identical to \texttt{saa}) or converges toward the solution of \texttt{ro} (when $\epsilon$ is sufficiently large to encompass the worst-case scenario, \eqref{prob} becomes identical to \texttt{ro} as the worst-case distribution will solely support the worst-case outcome). For intermediate values of $\epsilon$, \eqref{prob} tends to yield superior $x$ compared to \texttt{saa} and \texttt{ro}, especially in cases of limited sample size $n$. 

An optimal choice of $\epsilon$ for \eqref{prob} would be the one that excels in out-of-sample scenarios. With a sufficient number $M$ of testing scenarios, one can select $\epsilon$ by minimizing the expected system cost, approximated using the testing dataset denoted by $\{\zeta'_{ij}\}_{j=1}^{M}$. Formally, this involves finding $\epsilon$ that minimizes $c^T x(\epsilon) + \sum_{i\in \mc L}\sum_{j=1}^{M} \frac{1}{M} Z_i(x(\epsilon), \zeta'_{ij})$, where $x(\epsilon)$ represents the solution of \eqref{prob} given $\epsilon$. We refer to \eqref{prob} with the optimal choice of $\epsilon$ as \texttt{opt}.

However, as the true distribution $\mb Q_i$ is rarely known, and a testing dataset of sufficient size for evaluating out-of-sample performance is often unavailable, estimating the optimal $\epsilon$ relies mainly on the available training dataset. In this paper, we employ the holdout method for estimation, denoted by \texttt{hm}. The holdout method involves partitioning the training dataset into training and validation datasets, with $n_T$ and $n_V = n-n_T$ denoting their respective cardinalities. Utilizing only the new training dataset of size $n_T$ for each subnetwork $i$, we solve \eqref{prob} while varying $\epsilon$ to obtain $x(\epsilon)$. Subsequently, the total expected system cost of $x(\epsilon)$ is approximated using the validation dataset. Let $\epsilon'$ be the value that minimizes the estimated total system cost, and the corresponding $x(\epsilon')$ represents the solution of \texttt{hm}. In this paper, we allocate 20\% of the data for validation in the \texttt{hm} process. There exist various alternative methods for approximating the optimal $\epsilon$ developed in the statistical learning field, such as k-fold cross-validation used in \cite{esfahani2018data}.

}

\section{Test System Description}\label{appendix:test-system}
The test system has 8 line switches named $\{$sw1, sw2, $\cdots$, sw8$\}$ and 4 CBs at buses 83, 88, 90, and 92. The substation regulator is treated as an OLTC, in which the CC decides its voltage magnitude squared, $v_1$. Other regulators are removed, and the load transformers are modeled as lines with equivalent impedance. 
Since the original data is in three phases, we modify the data to get a single-phase test system as follows: (i) loads are averaged over phases, and (ii) the line impedance $r_{ij} + \i x_{ij}$ is obtained by the maximum diagonal element of its impedance matrix. Only the voltage source serves as a dispatchable generator. The voltage bounds are set to be [0.8 p.u., 1.2 p.u.].

{\color{black}\section{Figures}\label{sec:figs}
\begin{figure*}[t!]
    \centering
    \subfloat[$n=5$\label{a}]{\includegraphics[width=0.8\textwidth]{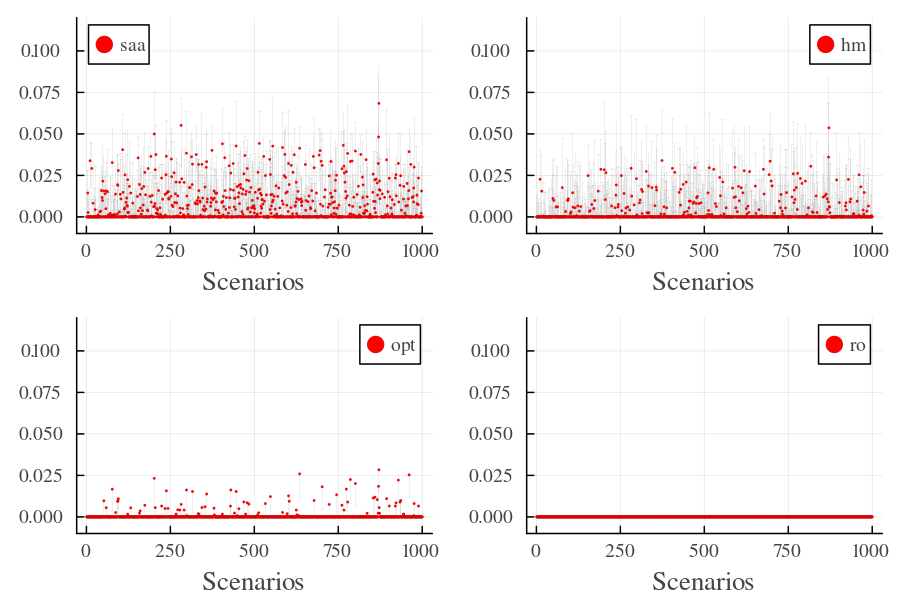}}
    \hfill
    \subfloat[$n=30$\label{b}]{\includegraphics[width=0.8\textwidth]{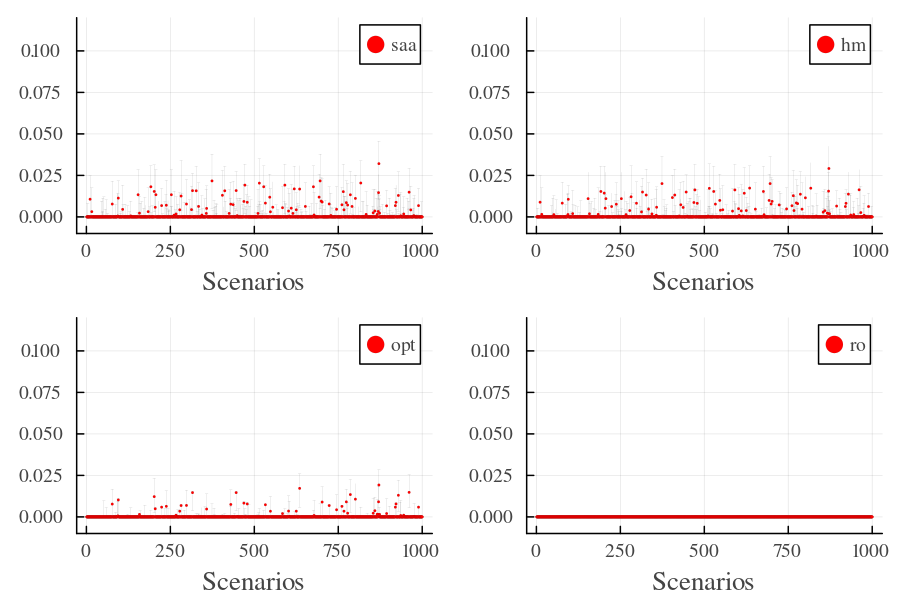}}
    \caption{Load shed in p.u. per scenario over 50 simulations}
\label{fig:loadshed_by_methods}
\end{figure*}

\begin{figure*}
    \centering
    \includegraphics[width=0.8\linewidth]{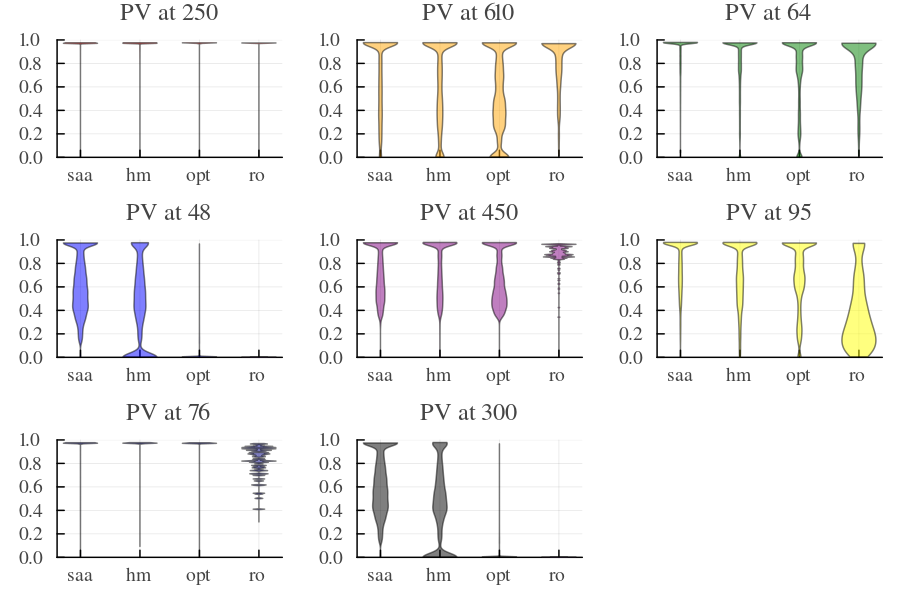}
    \caption{PV utilization, that is, $\sqrt{g_l^p+g_l^q}/\bar p_l \xi_l$ ($n=5$)}
    \label{fig:util_N5}
\end{figure*}
\begin{figure*}
    \centering
    \includegraphics[width=0.8\linewidth]{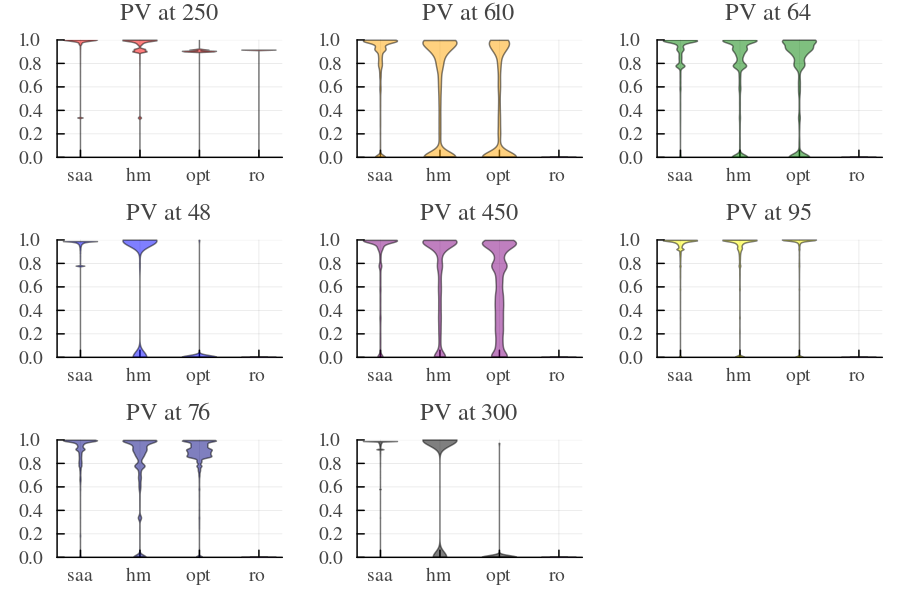}
    \caption{PV power factors, that is, $g_l^p/\sqrt{g_l^p+g_l^q}$ ($n=5$)}
    \label{fig:pf_N5}
    \vspace{-3mm}
\end{figure*}

}
\end{appendices} 

\end{document}